\documentclass{article}

\usepackage{microtype}
\usepackage{graphicx}
\usepackage{subfigure}
\usepackage{array}
\usepackage{arydshln}
\usepackage{lipsum}
\usepackage{booktabs} 
\usepackage{multirow, makecell}
\usepackage{amsmath}
\usepackage{amssymb}
\usepackage{amsthm}
\usepackage{amsfonts}
\usepackage{dsfont}
\usepackage{algorithmicx}
\usepackage{algorithm}
\usepackage{algpseudocode}

\usepackage{stackengine}

\usepackage{amsmath}
\DeclareMathOperator*{\argmax}{arg\,max}

\usepackage{hyperref}

\newtheorem{theorem}{Theorem}

\newtheorem{proposition}[theorem]{Proposition}
\newtheorem{corollary}{Corollary}

\graphicspath{{figures/}}
\usepackage[round, sort,comma,authoryear]{natbib}

\usepackage[]{authblk}

\usepackage{color}

\newcommand\rev[1]{\textcolor{black}{#1}}

\newcommand{\p}{\phantom{0}}
\newcommand{\pp}{\phantom{00}}

\newcommand{\pex}{\phantom{$^{(0)}$}}
\newcommand{\pexs}{\phantom{$^*$}}

\title{A model-free approach for solving choice-based competitive facility location problems using simulation and submodularity}

\author[1,2]{Robin Legault\thanks{Corresponding author}}
\author[2]{Emma Frejinger}

\affil[1]{Operations Research Center, Massachusetts Institute of Technology\\  \texttt{legault@mit.edu}\\
  }
\affil[2]{CIRRELT, Université de Montréal\\ \texttt{emma.frejinger@umontreal.ca}}
\date{}

\begin{document}

\maketitle

\begin{abstract}
This paper considers facility location problems in which a firm entering a market seeks to open facilities on a subset of candidate locations so as to maximize its expected market share, assuming that customers choose the available alternative that maximizes a random utility function. \rev{We introduce a deterministic equivalent reformulation of this stochastic problem as a maximum covering location problem with an exponential number of demand points, each of which is covered by a different set of candidate locations. Estimating the prevalence of these preference profiles through simulation generalizes a sample average approximation method from the literature and results in a maximum covering location problem of manageable size. To solve it, we develop a partial Benders reformulation in which the contribution to the objective of the least influential preference profiles is aggregated and bounded by submodular cuts. This set of profiles is selected by a knee detection method that seeks to identify the best tradeoff between the fraction of the demand that is retained in the master problem and the size of the model. We develop a theoretical analysis of our approach and show that the solution quality it provides for the original stochastic problem, its computational performance, and the automatic profile-retention strategy it exploits are directly connected to the entropy of the preference profiles in the population.} Computational experiments on existing and new benchmark sets indicate that our approach dominates the classical sample average approximation method on large instances of the competitive facility location problem, can outperform the best heuristic method from the literature under the multinomial logit model, and achieves state-of-the-art results under the mixed multinomial logit model. \rev{We characterize a broader class of problems, which includes assortment optimization, to which the solving methodology and the analyses developed in this paper can be extended.}\\

\end{abstract}

\section{Introduction}\label{sec:intro}
In a wide array of real-world management applications, the impact of a firm's decisions on the level of demand for its products or services depends on the preferences of a population of heterogeneous customers. The core assumption of utility maximization theory is that agents evaluate each available alternative and select the one that maximizes their utility function. Embedding random utility maximization (RUM) models into optimization problems leads to so-called choice-based optimization problems, prime examples of which include pricing \citep[e.g.,][]{davis2017pricing, gallego2014multiproduct, li2019product, paneque2022lagrangian} and assortment optimization problems \citep[e.g.,][]{liu2020assortment, rusmevichientong2010dynamic}.

The competitive facility location problem is another important problem that requires modeling customer demand at a disaggregate level. In the last two decades, a growing body of literature has proposed exact methods for solving the choice-based competitive facility location problem (CBCFLP) under the multinomial logit (MNL) model \citep{aros_2013, Benati_2002, Freire_2016, Ljubic_2018, Mai_2020, Zhang_2012}. A limitation of the MNL model is that it implies proportional substitution patterns, meaning that introducing a new alternative divides by the same factor the probability for each existing alternative to be selected by any given customer. This property leads to unrealistic demand representation in many cases. Although empirical studies underscore the need to capture unobserved taste variations and spatial correlation to faithfully model location-related behavior \citep{Bhat_2004, Miyamoto_2004, Muller_2012}, relatively little attention has been devoted to the study of the CBCFLP under less restrictive modeling assumptions. Notable exceptions include an exact branch-and-cut (B\&C) method for the nested logit choice rule \citep{mendez2023follower} and a heuristic local search approach that can be applied under any generalized extreme value (GEV) models \citep{dam2022submodularity}. In addition to the MNL and nested logit models \citep{Williams_1977}, prominent members of the GEV family include the paired combinatorial logit model \citep{Chu_1989} and the cross-nested logit model \citep{Vovsha_1997}. 

The most widely studied RUM model outside of the GEV family is the mixed multinomial logit (MMNL) model, which can approximate any RUM model with arbitrary precision \citep{McFadden_Train_2000}. Some facility location studies \citep[e.g.,][]{dam2022submodularity, Haase_2016, Mai_2020} have highlighted that this fully flexible model can be approximated by a MNL model through simulation. Unfortunately, obtaining near-optimal solutions for MMNL instances with this approach requires solving very large MNL instances in some cases, which can be computationally intractable even for state-of-the-art methods. 

Another approach that can accommodate flexible RUM models is the sample average approximation framework introduced by \citet{Haase_2009}. Initially developed for the MNL model and later extended to the MMNL model \citep{Haase_2013, Haase_2016}, this method approximates the utility function of each customer by sampling realizations of its random terms. \rev{The resulting approximation of the CBCFLP can be formulated as a maximum covering location problem (MCLP) \citep{church1974maximal} in which a demand point is generated for each simulated customer.} This method has been applied in school location planning \citep{haase2019facility} and to optimize electric vehicle charging station placement \citep{lamontagne2023optimising}. More generally, sample average approximation is receiving increasing attention in choice-based optimization as it makes it relatively straightforward to integrate advanced RUM models into mixed-integer linear programming (MILP) models \citep{Paneque_2021}. \rev{As opposed to model-specific methods that exploit the structure of a given family of RUM models, simulation-based approaches like the one we present in this paper do not make any restrictive assumption on the utility functions of the customers and are thus said to be model-free.}

\rev{Our first main contribution is to introduce a deterministic equivalent reformulation of the CBCFLP as a very large-scale MCLP in which each possible preference profile is represented by a unique demand point. The weight of these preference profiles can be approximated by simulation, and those with trivial estimated weight can be removed from the model, leading to a more compact yet equivalent reformulation of the sample average approximation model of \citet{Haase_2016}. Our second contribution is to propose a partial Benders decomposition \citep{crainic2021partial} of the MCLP in which the preference profiles with the highest weights are explicitly represented in the master problem, and the remaining ones are aggregated and replaced in the objective by an auxiliary variable that is bounded by submodular cuts. We propose to partition the profiles based on a knee detection method, leading to a profile-retention strategy that does not rely on user-defined parameters. Our third contribution is to propose an entropy measure to characterize the stochasticity level of an instance of the CBCFLP based on the level of concentration of the demand in the deterministic equivalent reformulation of the problem. We study the connections between this entropy measure, the quality of the sample average approximation of the original stochastic problem, and the performance of our partial Benders decomposition approach. Finally, we compare our approach with the sample average approximation of \citet{Haase_2016} and state-of-the-art model-specific exact \citep{Mai_2020} and heuristic \citep{dam2022submodularity} methods on existing MNL instances and new MMNL instances of the CBCFLP. Our computational study shows that our algorithm performs better than the classical sample average approximation method and the state-of-the-art heuristic method in terms of CPU time and solution quality for most MNL instances and achieves state-of-the-art results under the MMNL model.} 

The paper is structured as follows. Section~\ref{sec:max_capture_problem} presents the problem and introduces our notation. Section~\ref{sec:methods} discusses existing methods from the literature for solving the CBCFLP. Section~\ref{sec:simul_method} presents our \rev{partial Benders decomposition} approach. \rev{Section~\ref{sec:info_theo} presents the information-theoretic study.} The computational experiments are reported in Section~\ref{sec:experiments}, and Section~\ref{sec:conclusion} concludes the paper.

\section{Choice-Based Competitive Facility Location}\label{sec:max_capture_problem}
We consider a probabilistic facility location problem in a competitive market with utility-maximizing customers. Given a set $D$ of available locations on which a firm can install new facilities and a set $E$ of locations that competing facilities already occupy, the goal is to identify a feasible configuration $\boldsymbol{x} \in X$ that maximizes the expected market share captured by the new facilities. A binary variable $x_d$ indicates whether the company decides to open a facility at location $d \in D$. The feasible domain $X \subseteq \{0,1\}^{|D|}$ is specified by linear business constraints, such as the number of facilities that can be opened or the firm's budget. We denote by $C = D \cup E$ the set of both candidate and existing locations.

Customers select the open facility that maximizes their utility function. The attributes of the customers among the population, such as their location and personal preferences, are modeled by a random vector $\boldsymbol{\theta}$ with support $\Theta$. The impact of these attributes on the utility of each location is modulated by coefficients $\boldsymbol{\beta} \in B$. The specification of these coefficients depends on the model, and they are possibly random. In real applications, the corresponding values are estimated on data. Finally, a random term $\varepsilon_c$ affects the utility of each alternative $c \in C$. For a solution $\bar{\boldsymbol{x}} \in X$, the open alternatives are given by the set $C_{\bar{\boldsymbol{x}}} = D_{\bar{\boldsymbol{x}}} \cup E$, where $D_{\bar{\boldsymbol{x}}} = \{d \in D: \bar{x}_d=1\}$. The utility of each available alternative $c \in C$ is denoted by $u_c(\boldsymbol{\theta}, \boldsymbol{\beta}, \boldsymbol{\varepsilon})$. In the case of RUM models with additive error terms, this value decomposes as $u_c(\boldsymbol{\theta},\boldsymbol{\beta},\boldsymbol{\varepsilon}) = v_c(\boldsymbol{\theta},\boldsymbol{\beta}) + \varepsilon_c$. In general, the random vectors $\boldsymbol{\theta},\boldsymbol{\beta}$ and $\boldsymbol{\varepsilon}$ can be dependent, and $v_c: \Theta \times B \rightarrow \mathbb{R}$ can be any real-valued function for each location $c \in C$.

In most of the literature, it is assumed that the population is composed of a finite number of customers $n \in N$, each described by observed attributes $\boldsymbol{\theta}^n$. In this context, the distribution of $\boldsymbol{\theta}$ is given by the empirical distribution of the observations $\boldsymbol{\theta}^n$. In machine learning terminology, this interpretation underlies a \textit{conditional}, or \textit{discriminative} perspective on the attributes of the customers. In other words, the observed set of customers is considered fixed instead of being seen as a random sample drawn from an underlying probability distribution. For the sake of generality, the notation adopted in this paper instead proposes a \textit{generative} perspective, in which all of the customers' attributes are modeled as random variables. The problem can then be stated as maximizing the expected market share $F:2^D \rightarrow [0,1]$ based on the joint distribution of the random vectors $\boldsymbol{\theta}, \boldsymbol{\beta}$, and $\boldsymbol{\varepsilon}$. Assuming continuous probability distributions, it can be expressed as follows:
\begin{align}
     \label{model:exact}\max_{\substack{\boldsymbol{x} \in X}} F(D_{\boldsymbol{x}}) &= \max_{\substack{\boldsymbol{x} \in X}} \mathbb{P}_{\boldsymbol{\theta}, \boldsymbol{\beta}, \boldsymbol{\varepsilon}}\big[\argmax_{c \in C_{\boldsymbol{x}}}\{u_c(\boldsymbol{\theta},\boldsymbol{\beta},\boldsymbol{\varepsilon})\} \in D_{\boldsymbol{x}}\big], \\
    &= \max_{\substack{\boldsymbol{x} \in X}} \int_{\boldsymbol{\theta} \in \Theta}\mathbb{P}_{\boldsymbol{\beta},\boldsymbol{\varepsilon}| \boldsymbol{\theta}}\big[\argmax_{c \in C_{\boldsymbol{x}}}\{u_c(\boldsymbol{\theta},\boldsymbol{\beta},\boldsymbol{\varepsilon})\} \in D_{\boldsymbol{x}}\big] p(\boldsymbol{\theta})d\boldsymbol{\theta},\\
    &= \max_{\substack{\boldsymbol{x} \in X}} \int_{\boldsymbol{\theta} \in \Theta}\sum_{d \in D_{\boldsymbol{x}}}\mathbb{P}_{\boldsymbol{\beta},\boldsymbol{\varepsilon}| \boldsymbol{\theta}}\big[u_d(\boldsymbol{\theta},\boldsymbol{\beta},\boldsymbol{\varepsilon}) \geq u_{c}(\boldsymbol{\theta},\boldsymbol{\beta},\boldsymbol{\varepsilon}), \ \forall c \in C_{\boldsymbol{x}} \big] p(\boldsymbol{\theta})d\boldsymbol{\theta},
\end{align}
which, for RUM models with additive error terms, can be rewritten as:
\begin{align}
    \max_{\substack{\boldsymbol{x} \in X}} \int_{\boldsymbol{\theta} \in \Theta}\sum_{d \in D_{\boldsymbol{x}}}\mathbb{P}_{\boldsymbol{\beta},\boldsymbol{\varepsilon}| \boldsymbol{\theta}}\big[\varepsilon_{c}-\varepsilon_d \leq v_d(\boldsymbol{\theta},\boldsymbol{\beta})-v_{c}(\boldsymbol{\theta},\boldsymbol{\beta}), \ \forall c \in C_{\boldsymbol{x}} \big] p(\boldsymbol{\theta})d\boldsymbol{\theta}.  \label{eq:capProb}
\end{align}

For the choice probability in~(\ref{eq:capProb}) to be directly embedded in an optimization problem, the cumulative distribution function of the difference $\varepsilon_{c}-\varepsilon_d$ must admit an analytical form. This has motivated the development of increasingly flexible models belonging to the GEV family, which all result in closed-form expressions \citep{FosgMcFaBier13}. 
\rev{A} large part of the RUM literature assumes the coefficients $\boldsymbol{\beta}$ to be independent of the customers’ attributes $\boldsymbol{\theta}$ and is based on additive error terms $\varepsilon_c$ defined as independent and identically distributed (i.i.d.) type I extreme value (or Gumbel) random variables. These assumptions result in the MMNL model \rev{that has the nice property that it can approximate any RUM model (including the GEV family) with arbitrary precision \citep{McFadden_Train_2000}. The} probability of capturing the demand of a customer with attributes $\boldsymbol{\theta}$ is \rev{then} given by the following expression:
\begin{alignat}{1}
\label{eq:ProbMMNL}\mathbb{P}_{\boldsymbol{\beta},\boldsymbol{\varepsilon}| \boldsymbol{\theta}}\big[\varepsilon_c-\varepsilon_d \leq v_d(\boldsymbol{\theta},\boldsymbol{\beta})-v_c(\boldsymbol{\theta},\boldsymbol{\beta}), \ \forall c \in C_{\boldsymbol{x}} \big] = \int_{\boldsymbol{\beta} \in B} \frac{\sum_{d \in D_{\boldsymbol{x}}}e^{v_d(\boldsymbol{\theta},\boldsymbol{\beta})}}{\sum_{c \in C_{\boldsymbol{x}}}e^{v_c(\boldsymbol{\theta},\boldsymbol{\beta})}} p(\boldsymbol{\beta}) d\boldsymbol{\beta}.
\end{alignat}

Under the MMNL model and a generative perspective, model~(\ref{model:exact}) can be formulated as:
\begin{alignat}{2}
    \label{eq:MMNL} &\max_{\substack{\boldsymbol{x} \in X}} \int_{\boldsymbol{\theta} \in \Theta} \int_{\boldsymbol{\beta} \in B} \frac{\sum_{d \in D_{\boldsymbol{x}}}e^{v_d(\boldsymbol{\theta},\boldsymbol{\beta})}}{\sum_{c \in C_{\boldsymbol{x}}}e^{v_c(\boldsymbol{\theta},\boldsymbol{\beta})}} p(\boldsymbol{\beta})  p(\boldsymbol{\theta}) d\boldsymbol{\beta} d\boldsymbol{\theta}.
\end{alignat}

In turn, under the conditional perspective, it becomes: 
\begin{alignat}{2}
     \label{eq:MMNL_conditional}&\max_{\substack{\boldsymbol{x} \in X}} \sum_{n \in N} q_n \int_{\boldsymbol{\beta} \in B} \frac{\sum_{d \in D_{\boldsymbol{x}}}e^{v_d(\boldsymbol{\theta},\boldsymbol{\beta})}}{\sum_{c \in C_{\boldsymbol{x}}}e^{v_c(\boldsymbol{\theta},\boldsymbol{\beta})}} p(\boldsymbol{\beta}) d\boldsymbol{\beta},
\end{alignat}
where $q_n = p(\boldsymbol{\theta}^n)$ denotes the weight of each customer $n \in N$.

The MNL model is a specific case where the coefficients $\boldsymbol{\beta}$ are deterministic, which simplifies the generative MMNL model~(\ref{eq:MMNL}) to:
\begin{alignat}{2}
     \label{eq:MNL} &\max_{\substack{\boldsymbol{x} \in X}} \int_{\boldsymbol{\theta} \in \Theta} \frac{\sum_{d \in D_{\boldsymbol{x}}}e^{v_d(\boldsymbol{\theta},\boldsymbol{\beta})}}{\sum_{c \in C_{\boldsymbol{x}}}e^{v_c(\boldsymbol{\theta},\boldsymbol{\beta})}} p(\boldsymbol{\theta}) d\boldsymbol{\theta}.
\end{alignat}

Similarily, under the conditional perspective, (\ref{eq:MMNL_conditional}) simplifies to: 
\begin{alignat}{2}
     \label{eq:MNL_conditional}&\max_{\substack{\boldsymbol{x} \in X}} \sum_{n \in N} q_n \frac{\sum_{d \in D_{\boldsymbol{x}}}e^{v^n_{d}}}{\sum_{c \in C_{\boldsymbol{x}}}e^{v^n_{c}}},
\end{alignat}
where the perceived utility of an available alternative $c \in C$ for customer $n \in N$ is given by $v^n_{c}:=v_c(\boldsymbol{\theta}^n, \boldsymbol{\beta})$. In this case, and more generally if the support $\Theta$ is discrete and finite, the problem \rev{can} be solved exactly by nonlinear programming (MINLP) solvers or by model-specific algorithms \citep{Ljubic_2018, Mai_2020}. In order to solve MMNL or generative MNL instances approximately using these methods, the distribution of the random variables $\boldsymbol{\theta}$ and/or $\boldsymbol{\beta}$ can be approximated by the empirical distribution provided by a set of their realizations, resulting in a conditional MNL problem.

\section{Methods from the Literature}\label{sec:methods}
This section discusses the most recent state-of-the-art exact and heuristic methods for the CBCFLP under the MNL and MMNL models, as well as the sample average approximation method our approach builds on. 

\subsection{State-of-the-Art Exact Methods}\label{subsec:MOA}
\rev{There are currently two different algorithms that can be regarded as state-of-the-art exact methods for solving the CBCFLP under the MNL model and a conditional perspective. The first one, by \citet{Ljubic_2018}, is a B\&C procedure in which the contribution of each individual customer to the objective function is bounded separately by outer-approximation and submodular cuts. The second one, by \citet{Mai_2020}, is a cutting plane method that generates outer-approximation cuts for groups of customers at each iteration. The algorithm of \citet{Ljubic_2018} has shown to be the fastest of the two, by a factor of four on average, for the smallest benchmark instances of the literature. However, for the instances comprising the largest number of customers, the method of \citet{Mai_2020} has shown to be the fastest by about two orders of magnitude. This paper primarily focuses on instances with a large number of customers. In the following discussion, we thus present more details on the latter approach, to which our own method will be compared in the computational experiments.}

The multicut outer-approximation (MOA) algorithm \citep{Mai_2020} is based on the outer-approximation scheme \citep{Duran_1986}. In this method, model~(\ref{eq:MNL_conditional}) is first reformulated as a minimization problem, with objective function $G(x) = -\sum_{n \in N} q_n \frac{\sum_{d \in D_{\boldsymbol{x}}}e^{v^n_{d}}}{\sum_{c \in C_{\boldsymbol{x}}}e^{v^n_{c}}}$. The main idea of the MOA algorithm is to partition the set $N$ of customers into $T$ subsets so that $G(\boldsymbol{x})$ can be expressed as a sum of $T$ convex and continuously differentiable functions $g_t(\boldsymbol{x})$. At each iteration of the algorithm, the master problem $\min_{\boldsymbol{x} \in X}\{\sum_{t=1}^T \phi_t | \phi_t \geq L_t, \ \Pi_t\boldsymbol{x}-\boldsymbol{1}\phi_t \leq \pi_{0t} \ \forall t\}$ is solved, where the value of each component $g_t(\boldsymbol{x})$ is replaced by the decision variable $\phi_t$, $\Pi_t\boldsymbol{x}-\boldsymbol{1}\phi_t \leq \pi_{0t}$ is the set of subgradient cuts corresponding to $g_t(\boldsymbol{x})$, and $L_t$ is a lower bound on $g_t(\boldsymbol{x})$. Its optimal solution $\boldsymbol{x}^*$ is then used to add up to $T$ new subgradient cuts $\phi_t \geq \nabla g_t(\boldsymbol{x}^*)(\boldsymbol{x}-\boldsymbol{x}^*) + g_t(\boldsymbol{x}^*)$, $t=1,\dots,T$ to the master problem, where $\nabla g_t(\boldsymbol{x}^*)$ is the gradient of $g_t$ evaluated at $\boldsymbol{x}^*$. 

This multicut approach generalizes an earlier single-cut algorithm by \citet{Bonami_2008}. The single-cut version is obtained by using only $T=1$ objective function $g_1(\boldsymbol{x})=G(\boldsymbol{x})$, while taking $T \in \{2,\dots,|N|-1\}$ leads to what is sometimes referred to as a \textit{hybrid} approach \citep{Birge_1988} in the stochastic programming literature. Selecting a large value of $T$ tends to limit the number of iterations of MOA, but adding too many cuts at each iteration makes the master problem more challenging to solve. \citet{Mai_2020} report that neither the single-cut nor the pure multicut version with $T=|N|$ usually results in the most efficient version of the MOA. Furthermore, the best value of $T$ varies significantly across instances. Unfortunately, no efficient rule for automatically selecting $T$ is known, making the performance of MOA heavily dependent on a user-defined parameter.

\subsection{State-of-the-Art Heuristic Method}\label{subsec:GGX}
An efficient heuristic algorithm that can be applied under any GEV model has been proposed by \citet{dam2022submodularity}. This algorithm, called GGX (for Greedy heuristic, Gradient-based local search, and eXchanging), has been shown to identify an optimal solution to most MNL and MMNL instances from classical benchmark sets in a fraction of the time required by MOA to solve the problem to proven optimality. Encouraging computational results have also been observed on nested logit instances, although the performance of GGX on these instances is difficult to evaluate properly due to the lack of comparison with an exact method in this study. 

The GGX algorithm consists of three phases. First, starting from the trivial set $D_{\boldsymbol{x}}=\emptyset$, a greedy solution is constructed by repeatedly opening the location $d \in D \setminus D_{\boldsymbol{x}}$ leading to the largest increase in the objective value. Second, a local search within a region of increasing size based on a linear approximation of the objective function at the current solution is iteratively performed. The last phase is a greedy local search in which sets of closed and open locations are iteratively exchanged until a local maximum is found. \citet{dam2022submodularity} allow for at most two pairs of locations to be exchanged at each iteration.

Due to the monotonicity and submodularity of the CBCFLP under GEV models \citep{dam2022submodularity}, it follows from \citet{nemhauser1978analysis} that the greedy heuristic performed in the first phase of GGX is a $(1 - 1/e)$ approximation algorithm. In other words, the objective value of the solution returned by GGX is guaranteed to be at least $\approx 0.632$ times the optimal value. However, like most local search algorithms, GGX provides no stronger theoretical guarantee. \rev{In particular, the family of instances of the general covering problem presented in \citet{hochbaum1998analysis} for which the greedy algorithm reaches the worst-case approximation ratio can easily be adapted to the context of the CBCFLP to establish that the $(1 - 1/e)$ bound is also tight for GGX. Such instances are, however, very artificial and unlikely to appear in real applications. We illustrate a simpler and more realistic case in which GGX also leads to a strongly suboptimal solution.}

\begin{figure}[htbp]
\caption{Optimal solution (left, expected market share of 65.61\%) and solution obtained with GGX (right, expected market share of 58.66\%). The shade of a customer indicates the probability that they select a facility of the firm.}
\label{fig:GGX_counterexample}
\begin{center}
    \includegraphics[width=\textwidth]{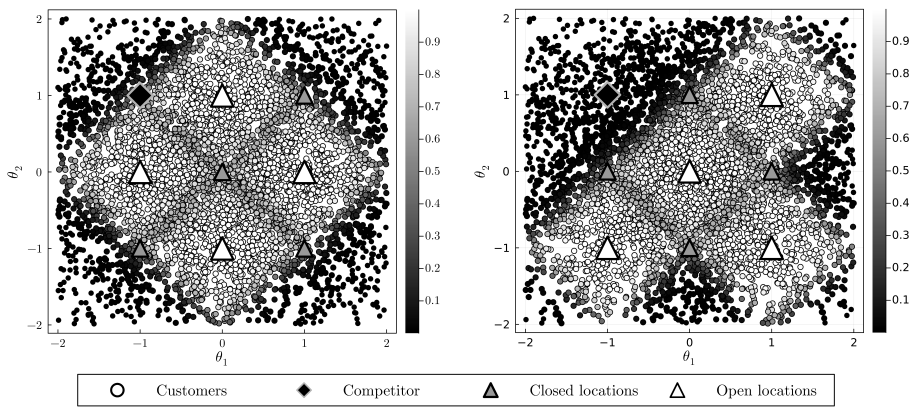}
\end{center}
\end{figure}

Figure \ref{fig:GGX_counterexample} illustrates the behavior of GGX on a conditional MNL instance with $|D|=8$ candidate locations, a budget of $4$ facilities, and $|E|=1$ existing location that provides the same perceived utility to each customer. The instance comprises $|N| = \text{5,000}$ customers generated according to a censored spherical normal centered at the origin. In this example, the central location is the first to be included in the greedy solution, followed by the three corner locations. This local maximum cannot be improved by the gradient-based local search of the second phase of GGX. Similarly, it cannot be improved by the greedy exchanging phase unless up to four pairs of locations can be exchanged at each iteration, which corresponds to performing an exhaustive search on the feasible domain $X$. The solution returned by GGX thus corresponds to opening all the locations that are closed in the optimal solution, and vice versa.

\subsection{Sample Average Approximation Method}\label{subsec:SAA}
In the sample average approximation approach introduced by \citet{Haase_2009} and studied by \citet{Haase_2013} and \citet{Haase_2016} in the context of the MMNL model under a conditional perspective, a set $S$ of scenarios $\{\boldsymbol{\beta}^s\}_{s \in S}$ is sampled from the distribution of the random coefficients $\boldsymbol{\beta}$. For each scenario $s \in S$ and each customer $n \in N$, a realization $\boldsymbol{\varepsilon}^{ns}$ of a vector of $|C|$ i.i.d. standard Gumbel variables is drawn. The problem associated with the empirical distribution $\{\boldsymbol{\theta}^n\}_{n \in N}$, $\{\boldsymbol{\beta}^s\}_{s \in S}$, $\{\boldsymbol{\varepsilon}^{ns}\}_{n \in N, s \in S}$ is given by:
\begin{alignat}{2}
     \label{model:simul_NS}&\max_{\substack{\boldsymbol{x} \in X}} \frac{1}{|S|}\sum_{n \in N} q_n \sum_{s \in S}\mathds{1}\big[\argmax_{c \in C_{\boldsymbol{x}}}\{u_c(\boldsymbol{\theta}^n,\boldsymbol{\beta}^s,\boldsymbol{\varepsilon}^{ns})\} \in D_{\boldsymbol{x}}\big].
\end{alignat}

Model~(\ref{model:simul_NS}) can be rewritten as the following 0-1 linear program:
\begin{alignat}{2}
     \label{model:SAA}\text{[SAA($N,S$)] }& &\max_{\substack{\boldsymbol{x} \in X \\ \boldsymbol{y} \in \{0,1\}^{|N|\times|S|} }} \hspace{0.3cm} \frac{1}{|S|} \sum_{n \in N} q_n \sum_{s \in S} y_{ns},\\
     \label{model:SAA:con}&\text{s.t. } &  y_{ns} \leq \sum_{d \in D} a^{ns}_{d}x_d, \hspace{0.5cm}&  \forall n \in N, \forall s \in S,
\end{alignat}
where the binary decision variable $y_{ns}$ indicates whether the simulated customer $(n,s) \in N \times S$ selects one of the new facilities of the firm. Binary coefficients $a^{ns}_{d} = \mathds{1}[u_d(\boldsymbol{\theta}^n, \boldsymbol{\beta}^s, \boldsymbol{\varepsilon}^{ns}) \geq u_{e}(\boldsymbol{\theta}^n, \boldsymbol{\beta}^s, \boldsymbol{\varepsilon}^{ns}) \ \forall e \in E]$ specify the candidate locations $d \in D$ that would be preferred by customer $n \in N$ to all the competing facilities $e \in E$ under scenario $s \in S$. Constraints~(\ref{model:SAA:con}) require the firm to open at least one such facility $d \in D$ to capture the demand of simulated customer $(n,s)$.

Since the coefficients $a^{ns}_{c}$ are computed outside of the optimization problem, SAA($N,S$)~(\ref{model:SAA})-(\ref{model:SAA:con}) is a 0-1 linear programming problem on a set of $|D| + |N|\cdot|S|$ decision variables. This formulation corresponds to a MCLP and is solved directly using a MILP solver in \citet{Haase_2016}.

As $|S| \rightarrow \infty$, model~(\ref{model:simul_NS}) approximates the conditional MMNL model~(\ref{eq:MMNL_conditional}) with an arbitrarily high precision. The asymptotic convergence property of sample average approximation provides probabilistic guarantees \citep{kim2015guide} that can make this approach preferable to greedy heuristics to avoid severely suboptimal solutions such as the one illustrated in Figure \ref{fig:GGX_counterexample}. However, since the number of decision variables and constraints of model~\text{SAA($N,S$)} grow linearly with $|S|$, the quality of the optimal solution's sample average estimate that can be obtained with a reasonable computational budget is often limited. 

\section{\rev{Simulation-Based Partial Benders Decomposition Method}}\label{sec:simul_method}
This section presents our model-free approach for the CBCFLP. Section~\ref{sec:simul_method:DEQ} introduces a deterministic equivalent reformulation of the problem. \rev{Section~\ref{sec:simul_method:simul} shows that approximating this model by simulation results in a compact MCLP. Section~\ref{sec:simul_method:PB} presents our partial Benders decomposition approach for solving this MCLP. Our method exploits a new profile-retention strategy, which is presented in Section~\ref{sec:simul_method:knee}. }

\subsection{Deterministic Equivalent Reformulation}\label{sec:simul_method:DEQ}
Our approach relies on a deterministic equivalent reformulation of model~(\ref{model:exact}):
\begin{alignat}{2}
     \label{model:DEQ1}\text{[DEQ] }& &\max_{\substack{\boldsymbol{x} \in X \\ \boldsymbol{y} \in \{0,1\}^{|P|} }} \hspace{0.3cm} \sum_{p \in P} \underline{\omega}_p y_{p},\\
     \label{model:DEQ1:con}&\text{s.t. } &  y_{p} \leq \sum_{d \in D} a^p_{d}x_d, \hspace{0.5cm}&  \forall p \in P.
\end{alignat}
\rev{In model DEQ, a demand point is associated with each \textit{preference profile} customers may have with positive probability. In our context, the preference profile of a customer corresponds to the subset of candidate locations they would prefer to the competing facilities. Each preference profile $p \in P$ can thus be encoded by the binary vector $\boldsymbol{a}^p \in \{0,1\}^{|D|}$ whose element $d \in D$ is set to 1 if the candidate location $d \in D$ covers the demand of $p$, and to 0 otherwise. The binary decision variable $y_p$ indicates whether at least one location $d$ such that $a^p_d = 1$ is opened in configuration $\boldsymbol{x} \in X$. Let the function $a_d: \Theta \times B \times \mathds{R} \rightarrow \{0,1\}$, defined as $a_d(\boldsymbol{\theta}, \boldsymbol{\beta}, \boldsymbol{\varepsilon}) = \mathds{1}[u_d(\boldsymbol{\theta}, \boldsymbol{\beta}, \boldsymbol{\varepsilon}) \geq u_{e}(\boldsymbol{\theta}, \boldsymbol{\beta}, \boldsymbol{\varepsilon}) \ \forall e \in E]$, indicate whether a customer defined by the vectors $\boldsymbol{\theta}, \boldsymbol{\beta}$ and $\boldsymbol{\epsilon}$ would patronize location $d \in D$ over the competitors, i.e., the opt-out alternatives.} The weight $\underline{\omega}_{p}$ is the probability $\mathbb{P}_{\boldsymbol{\theta},\boldsymbol{\beta},\boldsymbol{\varepsilon}}\left[ a_d(\boldsymbol{\theta}, \boldsymbol{\beta}, \boldsymbol{\varepsilon}) = a^{p}_d \ \forall d \in D \right]$ for a randomly selected customer to exhibit preference profile $\boldsymbol{a}^p$. Normalizing these weights lead to the unit vector $\boldsymbol{\omega}$, whose element $p \in P$ is given by $\omega_p = \mathbb{P}_{\boldsymbol{\theta},\boldsymbol{\beta},\boldsymbol{\varepsilon}}\left[ a_d(\boldsymbol{\theta}, \boldsymbol{\beta}, \boldsymbol{\varepsilon}) = a^{p}_d \ \forall d \in D | a_d(\boldsymbol{\theta}, \boldsymbol{\beta}, \boldsymbol{\varepsilon}) \neq \boldsymbol{0} \right] = \underline{\omega}_p/\sum_{p' \in P}\underline{\omega}_{p'}$. The following proposition shows the validity of this reformulation.
\begin{proposition}
\label{prop:DEQ_valid}
DEQ is a valid reformulation of model~(\ref{model:exact}).
\end{proposition}
\proof
The set $X$ of feasible configurations is the same for both models. Therefore, it suffices to show that the optimal value of DEQ given that the decision vector $\boldsymbol{x}$ is fixed to a feasible configuration $\bar{\boldsymbol{x}} \in X$ is equal to the objective value $F(D_{\bar{\boldsymbol{x}}})$ of the original model:
\begingroup
\allowdisplaybreaks
\begin{align}
\label{model:DEQ_restricted} &\max_{\substack{\boldsymbol{y} \in \{0,1\}^{|P|} }}\left\{\sum_{p \in P} \underline{\omega}_p y_p \Bigg| y_{p} \leq \sum_{d \in D} a^p_{d}\bar{x}_d, \ \forall p \in P\right\}, \\
    \label{model:DEQ_restricted_opt}=&\sum_{p \in P} \underline{\omega}_p \mathds{1}\left[ \sum_{d \in D} a^p_{d}\bar{x}_d \geq 1 \right], \\
    \label{eq:DEQ_restricted_def_omega} =&\sum_{p \in P} \mathbb{P}_{\boldsymbol{\theta},\boldsymbol{\beta},\boldsymbol{\varepsilon}}\left[ a_d(\boldsymbol{\theta}, \boldsymbol{\beta}, \boldsymbol{\varepsilon}) = a^{p}_d \ \forall d \in D \right] \mathds{1}\left[ \sum_{d \in D} a^p_{d}\bar{x}_d \geq 1 \right], \\
    \label{eq:DEQ_restricted_def_omega_before_tot_prob} =& \mathbb{P}_{\boldsymbol{\theta}, \boldsymbol{\beta}, \boldsymbol{\varepsilon}}\big[\sum_{d \in D}a_d(\boldsymbol{\theta}, \boldsymbol{\beta}, \boldsymbol{\epsilon})\bar{x}_d \geq 1 \big],\\
    =& \mathbb{P}_{\boldsymbol{\theta}, \boldsymbol{\beta}, \boldsymbol{\varepsilon}}\big[\sum_{d \in D_{\bar{\boldsymbol{x}}}}a_d(\boldsymbol{\theta}, \boldsymbol{\beta}, \boldsymbol{\epsilon}) \geq 1 \big],\\
    =& \mathbb{P}_{\boldsymbol{\theta}, \boldsymbol{\beta}, \boldsymbol{\varepsilon}}\big[\argmax_{c \in C_{\bar{\boldsymbol{x}}}}\{u_c(\boldsymbol{\theta},\boldsymbol{\beta},\boldsymbol{\varepsilon})\} \in D_{\bar{\boldsymbol{x}}}\big],\\
    =& F(D_{\bar{\boldsymbol{x}}}).
\end{align}
\endgroup
The restricted model DEQ with $\boldsymbol{x}=\bar{\boldsymbol{x}}$ is given by~(\ref{model:DEQ_restricted}) and is maximized by setting each decision variable $y_p$ to its maximum feasible value, hence equation~(\ref{model:DEQ_restricted_opt}). We obtain (\ref{eq:DEQ_restricted_def_omega}) by replacing the weights $\underline{\omega}_p$ with their definition. Expression (\ref{eq:DEQ_restricted_def_omega}) can be derived from (\ref{eq:DEQ_restricted_def_omega_before_tot_prob}) by applying the law of total probability. The last three equations are obtained by subsequently applying the definition of set $D_{\bar{\boldsymbol{x}}}$, coefficients $a_d(\boldsymbol{\theta}, \boldsymbol{\beta}, \boldsymbol{\epsilon})$ and function $F(D_{\bar{\boldsymbol{x}}})$.
\endproof

\subsection{\rev{Sample Average Approximation}}\label{sec:simul_method:simul}
In general, computing exactly the vector of weights $\underline{\boldsymbol{\omega}}$ \rev{in model DEQ} is not possible, as it requires to evaluate the following expression for each profile $p \in P$:
\begin{align}
    \underline{\omega}_p &= \int\limits_{\boldsymbol{\theta} \in \Theta} \int\limits_{\boldsymbol{\beta} \in B} \int\limits_{\boldsymbol{\varepsilon} \in \mathds{R}^{|C|}} \mathds{1}\left[ a_d(\boldsymbol{\theta}, \boldsymbol{\beta}, \boldsymbol{\varepsilon}) = a^p_{d} \ \forall d \in D\right]  p(\boldsymbol{\theta}, \boldsymbol{\beta}, \boldsymbol{\varepsilon}) d\boldsymbol{\varepsilon} d\boldsymbol{\beta} d\boldsymbol{\theta}.
\end{align}
However, as long as it is possible to draw samples from the joint distribution of the random variables $(\boldsymbol{\theta}, \boldsymbol{\beta}, \boldsymbol{\varepsilon})$, $\underline{\boldsymbol{\omega}}$ can be estimated by simulation. For a set $R$ of realizations $\{(\boldsymbol{\theta}^r, \boldsymbol{\beta}^r, \boldsymbol{\varepsilon}^r)\}_{r \in R}$, the estimated coefficients $\underline{\omega}_p$, $p \in P$, are given by:
\begin{align}
    \label{eq:omega_hat} \hat{\underline{\omega}}_p &= \frac{1}{|R|} \sum_{r \in R} \mathds{1}\left[ a_d(\boldsymbol{\theta}^r, \boldsymbol{\beta}^r, \boldsymbol{\varepsilon}^r) = a^p_{d} \ \forall d \in D\right],
\end{align}
which leads to the following approximation of DEQ:
\begin{alignat}{2}
     \label{model:DEQ_hat}\text{[$\widehat{\text{DEQ}}$]}& &\max_{\substack{\boldsymbol{x} \in X \\ \boldsymbol{y} \in \{0,1\}^{|\hat{P}|} }} \hspace{0.3cm} \sum_{p \in \hat{P}} \hat{\underline{\omega}}_p y_{p},\\
     \label{model:DEQ_hat:con}&\text{s.t. } &  y_{p} \leq \sum_{d \in D} a^p_{d}x_d, \hspace{0.5cm}&  \forall p \in \hat{P}.
\end{alignat}
Here, $\hat{P} \subseteq P$ is the set of profiles that have been observed at least once over the sample of simulated customers $R$, i.e., such that $\hat{\underline{\omega}}_p > 0$.

This simulation framework generalizes the classical sample average approximation approach. Indeed, the simulation method of \citet{Haase_2016} corresponds to building the set of realizations $R = N \times S$, with $\{(\boldsymbol{\theta}^r, \boldsymbol{\beta}^r, \boldsymbol{\varepsilon}^r)\}_{r \in R} = \{(\boldsymbol{\theta}^n, \boldsymbol{\beta}^s, \boldsymbol{\varepsilon}^{ns})\}_{(n,s) \in N \times S}$. Model SAA($N,S$) can then be reformulated as $\widehat{\text{DEQ}}$ by computing the weight of each profile $p \in \hat{P}=\{p \in P: \exists (n,s)\in N \times S \text{ such that } \boldsymbol{a}^p = \boldsymbol{a}^{ns} \}$ as:
\begin{equation}
    \label{eq:SAA->DEQ} \hat{\underline{\omega}}_p = \sum_{\substack{(n,s) \in N \times S \\ \boldsymbol{a}^{ns} = \boldsymbol{a}^p}} \frac{q_n}{|S|}.
\end{equation}
The following proposition shows that model $\widehat{\text{DEQ}}$ generalizes SAA($N,S$). 
\begin{proposition}
\label{prop:DEQ_gen_SAA}
The optimal value and the optimal configuration $\boldsymbol{x}^* \in X$ of model~$\widehat{\text{DEQ}}$, where the weights $\hat{\underline{\omega}}_p$ are given by (\ref{eq:SAA->DEQ}), are identical to those of SAA($N,S$).
\end{proposition}
\proof
The proof is similar to that of Proposition \ref{prop:DEQ_valid}. It is thus deferred to Online Appendix~A.
\endproof

The resulting $\widehat{\text{DEQ}}$ model contains fewer decision variables and constraints than SAA($N,S$) as it aggregates under a single profile $p \in P$ all the simulated customers $(n,s) \in N \times S$ sharing the same preferences $\boldsymbol{a}^{ns}=\boldsymbol{a}^p$. We call SAAA (for Sample Average Approximation with Aggregation) the method that consists in solving $\widehat{\text{DEQ}}$ directly using a solver, and denote by SAA the sample average approximation method of \cite{Haase_2016}. 

\subsection{\rev{Partial Benders Decomposition}}\label{sec:simul_method:PB}
\rev{In the context of the CBCFLP}, the number of possible preference profiles grows exponentially with the number of candidate locations. Solving $\widehat{\text{DEQ}}$ directly can thus be expected to become inefficient when $|D|$ is large and a high number of realizations $|R|$ is observed. In particular, for RUM models that include independent unbounded additive error terms, such as the MNL model, each of the $2^{|D|}-1$ non-trivial preference profiles can occur with positive probability, even in a conditional setting with only one customer. Most profiles $p \in \hat{P}$, however, typically have a negligible weight $\hat{\underline{\omega}}_{p}$ and thus only a small impact on the objective function. \rev{This property is exploited by the following reformulation, which aims at reducing} the computational sensitivity of our method to the number of observed preference profiles by aggregating the least important ones into a single composite customer. 

In model $\widehat{\text{DEQ}}(\hat{P}_1)$, only a subset $\hat{P}_1$ of the preference profiles $p \in \hat{P}$ are explicitly represented in the objective function. The captured market share on the remaining profiles $p \in \hat{P}_2=\hat{P} \setminus \hat{P}_1$ is represented by an auxiliary decision variable $\nu$\rev{:}
\begin{align}
    \label{model:DEQ1P1}\text{[$\widehat{\text{DEQ}}(\hat{P}_1)$] } &\max_{\substack{\boldsymbol{x} \in X \\ \boldsymbol{y} \in \{0,1\}^{|\hat{P}_1|} }} \hspace{0.3cm} \sum_{p \in \hat{P}_1} \hat{\underline{\omega}}_p y_{p} + \nu\\
     \label{model:DEQ1P1:con}\text{s.t. } &  y_{p} \leq \sum_{d \in D} a^p_{d}x_d, \hspace{0.5cm}&  \forall p \in \hat{P}_1,\\
     \label{model:DEQ1P1:SM1}& \nu \leq f\left(D_{\bar{\boldsymbol{x}}}\right) + \sum_{d \notin D_{\bar{\boldsymbol{x}}}}\rho_d\left(D_{\bar{\boldsymbol{x}}}\right)x_d - \sum_{d \in D_{\bar{\boldsymbol{x}}}}\rho_d(D \setminus \{d\})(1-x_d), & \forall \bar{\boldsymbol{x}} \in X.
\end{align}

\rev{For a set of open locations $D_{\bar{\boldsymbol{x}}} \subseteq D$, the contribution to the objective value of the set of profiles $\hat{P}_2$ is given by $f\left(D_{\bar{\boldsymbol{x}}}\right) = \sum_{p \in \hat{P}_2}\hat{\underline{\omega}}_p\mathds{1}\left[\sum_{d \in D{\bar{\boldsymbol{x}}}} a^p_{d} \geq 1\right]$. The function $f:2^{D}\rightarrow \mathds{R}_{\geq 0}$,  corresponds to the weighted set-union operator \citep{schrijver2003combinatorial} with positive weights $\hat{\underline{\omega}}_p$, and is thus submodular. It follows \citep[see][]{nemhauser1981maximizing, Ljubic_2018} that bounding $\nu$ by the submodular cuts (\ref{model:DEQ1P1:SM1}), where $\rho_d(D_{\bar{\boldsymbol{x}}}) = f\left(D_{\bar{\boldsymbol{x}}} \cup \{d\}\right) -f\left(D_{\bar{\boldsymbol{x}}}\right)$ denotes the marginal gain of opening a new location $d \in D \setminus D_{\bar{\boldsymbol{x}}}$, leads to a valid reformulation of model $\widehat{\text{DEQ}}$. In addition, since we have shown in Proposition \ref{prop:DEQ_valid} that the CBCFLP can equivalently be reformulated as a MLCP, it follows that the objective function of the CBCFLP is submodular under any RUM model. This observation generalizes the results of \citet{benati1997submodularity} and \citet{dam2022submodularity}, who proved the submodularity of the problem's objective function under the MNL model and GEV models, respectively. No proof of submodularity that directly applies to other families of RUM models, such as the MMNL model, was previously available in the CBCFLP literature. We note, however, that the same conclusion can also be reached as a consequence of the more general submodularity property of coverage-type objective functions under RUM models developed by \cite{berbeglia2020assortment} in the context of assortment optimization.}

\rev{It has been shown by \citet{coniglio2022submodular} that the submodular cuts (\ref{model:DEQ1P1:SM1}) are a special family of Benders cuts introduced by \citet{cordeau2019benders} for the MCLP. Consequently, model $\widehat{\text{DEQ}}(\emptyset)$ corresponds to a Benders decomposition reformulation of $\widehat{\text{DEQ}}$. Setting $\hat{P}_1 \notin \{\emptyset, \hat{P}\}$ instead leads to a partial Benders decomposition reformulation. This methodology has been extensively studied by \citet{crainic2021partial} in the context of two-stage stochastic multicommodity network design models. Its key idea is to include explicit information from the scenario subproblems in the master problem to reduce the number of feasibility cuts and improve the bound provided by the optimality cuts. Model $\widehat{\text{DEQ}}$ can be regarded as a two-stage stochastic problem where each preference profile constitutes a scenario. Selecting the set $\hat{P}_1$ to include in the master problem thus corresponds to what \citet{crainic2021partial} refer to as a scenario-retention strategy. In our case, we refer to it as a profile-retention strategy, and we introduce it in the following section.}

\rev{The following family of submodular cuts \citep{nemhauser1981maximizing} can replace the cuts (\ref{model:DEQ1P1:SM1}) or be added to model $\widehat{\text{DEQ}}(\hat{P}_1)$ to obtain an alternative valid reformulation of model $\widehat{\text{DEQ}}$:}
\begin{equation}
    \label{model:DEQ1P1:SM2} \nu \leq f\left(D_{\bar{\boldsymbol{x}}}\right) + \sum_{d \notin D_{\bar{\boldsymbol{x}}}}\rho_d\left(\emptyset\right)x_d - \sum_{d \in D_{\bar{\boldsymbol{x}}}}\rho_d(D_{\bar{\boldsymbol{x}}} \setminus \{d\})(1-x_d), \hspace{0.5cm} \forall \bar{\boldsymbol{x}} \in X.
\end{equation}
\rev{In addition, \citet{coniglio2022submodular} have shown that the submodular cuts (\ref{model:DEQ1P1:SM2}) can be strengthened and made equivalent to another family of Benders cuts which, for model $\widehat{\text{DEQ}}(\hat{P}_1)$, can be expressed as:}
\begin{equation}
    \label{model:DEQ1P1:BC2} \rev{\nu \leq f\left(D_{\bar{\boldsymbol{x}}}\right) + \sum_{d \notin D_{\bar{\boldsymbol{x}}}}\left(\sum_{p \in \hat{P}_2}\hat{\underline{\omega}}_pa^p_{d}\mathds{1}\left[\sum_{d \in D{\bar{\boldsymbol{x}}}} a^p_{d} \leq 1\right]\right) x_d - \sum_{d \in D_{\bar{\boldsymbol{x}}}}\rho_d(D_{\bar{\boldsymbol{x}}} \setminus \{d\})(1-x_d), \hspace{0.5cm} \forall \bar{\boldsymbol{x}} \in X.}
\end{equation}

\rev{In the strengthened cut (\ref{model:DEQ1P1:BC2}) associated with a solution $\bar{\boldsymbol{x}} \in X$, the customers that are covered by two or more locations $d \in D_{\bar{\boldsymbol{x}}}$ are removed from the set $\hat{P}_2$ on which the marginal gain $\rho_d(\emptyset)$ of opening a single location $d \notin D_{\bar{\boldsymbol{x}}}$ is computed in the first summation of (\ref{model:DEQ1P1:SM2}). New coefficients for the variables $x_d, d \notin D_{\bar{\boldsymbol{x}}}$ thus have to be computed each time a new cut is generated in (\ref{model:DEQ1P1:BC2}), whereas the values $\rho_d\left(\emptyset\right)$ only have to be computed once in the case of the submodular cuts (\ref{model:DEQ1P1:SM2}). In our numerical results, we observed that using cuts (\ref{model:DEQ1P1:BC2}) instead of (\ref{model:DEQ1P1:SM2}) substantially reduces the number of cuts that are generated in the solving process for MNL instances in which customers' preferences have a high level of stochasticity. However, the additional cost of computing the strengthened cuts exceeds the advantage provided by the improved bounds in terms of overall computing time for the largest instances of our testbed. The computational results presented in the paper are therefore based on the submodular cuts~(\ref{model:DEQ1P1:SM1}) and~(\ref{model:DEQ1P1:SM2}).}

Both sets of cuts are implemented using the \texttt{lazy-cut callback} routine of CPLEX and are applied globally in the branch-and-bound tree each time an integer solution violating them is found. Constraint $\nu \leq \sum_{p \in \hat{P}_2}\hat{\underline{\omega}}_p$ is added to the model to provide an initial valid upper bound on $\nu$.

\subsection{\rev{Profile-Retention Strategy}}\label{sec:simul_method:knee}
\rev{Our partial Benders decomposition} approach is agnostic to the partition of $\hat{P}$ regarding solution quality, as solving $\widehat{\text{DEQ}}(\hat{P}_1)$ yields an optimal solution to $\widehat{\text{DEQ}}$ for any $\hat{P}_1 \subseteq P$. However, its computational efficiency is highly dependent on this choice. \rev{The classical \textit{row covering} and \textit{convex hull} scenario-retention strategies developed by \citet{crainic2021partial} for partial Benders decomposition of two-stage stochastic programs primarily aim at reducing the number of feasibility cuts that have to be generated in the solving process. However, such considerations do not apply in our case because of the complete recourse structure of the MCLP. Instead, we introduce in this section a new strategy whose goal is to limit the number of decision variables of the model while minimizing the number of optimality cuts -- in our case, the submodular cuts (\ref{model:DEQ1P1:SM1}) and (\ref{model:DEQ1P1:SM2}) -- that have to be generated in the B\&C tree.}

Whereas a large part of the observed preference profiles must be included in $\hat{P}_2$ to significantly reduce the number of decision variables, a high number of submodular cuts may have to be generated in the B\&C tree if the auxiliary variable $\nu$ aggregates an excessively large fraction of the demand. This motivates the \rev{retention in the master problem} of the profiles that contribute the most to the objective function. Indexing the observed preference profiles by $p_1, \dots, p_{|\hat{P}|}$, where $\hat{\underline{\omega}}_{p_1} \geq \hat{\underline{\omega}}_{p_2} \geq \dots \geq \hat{\underline{\omega}}_{p_{|\hat{P}|}}$, we thus set $\hat{P}_1 = \{p_1, \dots, p_{i^*}\}$ for a given number of profiles $i^* \in \{0,\dots,|\hat{P}|\}$.

As $i^*$ increases, including additional profiles in $\hat{P}_1$ provides a decreasing marginal gain in the ratio $\Omega = \sum_{p \in \hat{P}_1}\hat{\underline{\omega}}_p / \sum_{p \in \hat{P}}\hat{\underline{\omega}}_p$ of the demand that is explicitly represented in the \rev{master problem}. Selecting an adequate cardinality for $\hat{P}_1$ thus corresponds to fixing an appropriate cutoff point in an increasing function with diminishing returns. A prime example of this type of problem arises in clustering tasks, where the marginal increase in the explained variation of the data decreases with the number of clusters. In this area, the number of clusters is usually determined based on a \textit{knee detection method}, where the knee of a function is defined as the maximizer of a curvature measure \citep{salvador2004determining}. 

\rev{The standard definition of curvature as the reciprocal of the radius of curvature of a space curve only applies to continuous functions \citep[see, for example,][]{do2016differential}.} However, a simple approach for approximating the point of maximum curvature for discrete data sets has been proposed by \citet{satopaa2011finding}. After normalizing the points in the unit square, the so-called \textit{Kneedle} method defines the knee as the point whose distance between the $y$-axis and $x$-axis coordinates is maximal. Formally, for a set of points $\{(i, j_i), i \in \{0,\dots,n\}\}$ respecting $j_i < j_{i+1} \ \forall i \in \{0,\dots,n-1\}$ and $(j_{i+2} - j_{i+1}) \leq (j_{i+1} - j_i) \ \forall i \in \{0,\dots,n-2\}$, Kneedle selects a point $(i^*, j_{i^*})$ such that:
\begin{equation} \label{eq:knee}
    i^* \in I^* = \argmax_{i \in \{0,\dots,n\}} \frac{j_i- j_0}{j_n - j_0} - \frac{i}{n}.
\end{equation}

This approach can be applied directly to our parameter selection problem by considering the points $\{(i, \Omega_i), i \in \{0,\dots,|\hat{P}|\}\}$, where $\Omega_i = \sum_{k = 1}^i \hat{\omega}_{p_k}$. It leads to selecting a point that maximizes the difference $\delta_i = \Omega_i - i/|\hat{P}|$ \rev{between the proportion of the demand $\Omega_i$ retained in the master problem and the proportion $i/|\hat{P}|$ of the original decision variables $y_p$ representing this demand.} When $|I^*|>1$, we select $i^* = \max I^*$. We observed that \rev{this profile-retention strategy} consistently provides a good approximation of the optimal cardinality of $\hat{P}_1$ for instances that can benefit from the \rev{partial Benders decomposition}. The characterization of such instances is discussed in Section~\ref{sec:info_theo}.

\rev{Selecting the profiles to retain in the master problem based on the knee detection method} also has the practical advantage of being independent of any user-defined parameter. This makes the performance comparison of our approach with other methods more objective. Henceforth, we denote the optimal value of the knee detection problem and the resulting relative weight of $\hat{P}_1$ by $\delta^*:=\delta_{i^*}$ and $\Omega^*:=\Omega_{i^*}$, respectively. \rev{We refer to the version of the partial Benders decomposition that is based on the Kneedle profile-retention strategy as PBD.} An example of the knee detection problem is presented in Online Appendix~B.

\section{Information-Theoretic Analysis}\label{sec:info_theo}
\rev{In this section, we adopt an information-theoretic perspective to study the discrepancy between the original CBCFLP model and its sample average approximation, as well as the performance of our method. Section~\ref{sec:info_theo:entropy_measure} presents the entropy measure on which our analysis is based. Sections~\ref{sec:info_theo:quality}~and~\ref{sec:info_theo:performance} study the impact of entropy on the solution quality provided by
simulation-based methods and on the computational
performance of PBD, respectively. In Section~\ref{sec:info_theo:discussion}, we discuss a class of problems to which our information-theoretic analysis and the general methodology developed in this paper can directly be extended.}

\subsection{\rev{Entropy Measure}}\label{sec:info_theo:entropy_measure}
\rev{In the deterministic equivalent reformulation DEQ of the CBCFLP, the demand is represented by a finite set $p \in P$ of preference profiles with normalized weights $\omega_p$. The distribution of these weights defines a categorical random variable $W$ with support $P$ and probability mass function (p.m.f.) $\boldsymbol{\omega}$. Its entropy is given by $H(W) = -\sum_{p \in P}\omega_p \log(\omega_p)$. In the same way, a categorical random variable $\hat{W}$ can be derived from the normalized weights $\hat{\omega}_p$ that indicate the importance of each observed preference profile $p \in \hat{P}$ in the sample average approximation model $\widehat{\text{DEQ}}$. Its entropy is given by $H(\hat{W}) = -\sum_{p \in \hat{P}}\hat{\omega}_p \log(\hat{\omega}_p)$.} The entropy of these random variables is an indicator of the level of concentration of the demand across the possible profiles $P$ (in the case of model DEQ) and across the observed profiles $\hat{P}$ (in the case of model $\widehat{\text{DEQ}}$). An instance of the CBCFLP achieves the lowest demand concentration when all the preference profiles $p \in P$ share the same normalized weight $\omega_p = \frac{1}{|P|}$. In this case, $W$ follows the discrete uniform distribution over $P$, and $H(W)=\log|P|$ reaches the theoretical upper bound on the entropy of a discrete random variable defined on a finite support of cardinality $|P|$. On the contrary, if a unique profile $p' \in P$ concentrates practically all the demand, say $\omega_{p'} = 1 - \epsilon$, then $H(W) \xrightarrow{\epsilon \rightarrow 0} 0$.

\rev{In the context of choice-based optimization, the definition of the utility function customers maximize is reflected in the entropy of the preference profiles. Whereas a customer who randomly ranks the opt-out alternative and the candidate locations exhibits each of the $|P|$ possible preference profiles with the same probability, a deterministic customer behaves according to a unique preference profile. Defining instances of the CBCFLP based on a single of these customers would lead to the two extreme cases previously discussed. In general, for a fixed set of customers, the entropy increases with the stochasticity of the RUM model. The entropy is thus a simple measure that can be used to characterize the structure of the demand. This leads to two main considerations.}

First, \rev{in a high entropy regime where} all the possible preference profiles are nearly equiprobable, then the coefficients of model $\widehat{\text{DEQ}}$ should converge rather slowly to those of DEQ. In such cases, we can expect simulation-based methods to produce suboptimal solutions for the original problem unless a very large sample of preference profiles is drawn.

Second, solving $\widehat{\text{DEQ}}$ to proven optimality should be easier if there exists a small set of preference profiles that covers a large proportion of the demand, as this would allow pruning earlier in the solving process the solutions that do not capture the most important profiles. \rev{PBD} explicitly exploits this structure through its \rev{profile-retention strategy} by including the most influential profiles in $\hat{P}_1$ and correcting the objective function for the remaining profiles $\hat{P} \setminus \hat{P}_1$ using submodular cuts. \rev{PBD} is thus expected to offer the most significant computational gain over SAAA when the demand is sufficiently concentrated.

\subsection{Impact of Entropy on Solution Quality}\label{sec:info_theo:quality}
The implicit assumption underlying any simulation approach in choice-based optimization is that the stochastic model describing the population's choices can be approximated efficiently from a sample of observed or reported preferences. In the case of model~(\ref{model:exact}), it means that it should be possible to faithfully approximate the coefficients $\underline{\omega}_p$, $p \in P$ of the deterministic equivalent formulation DEQ based on a finite set of realizations $\{(\boldsymbol{\theta}^r, \boldsymbol{\beta}^r, \boldsymbol{\varepsilon}^r)\}_{r \in R}$ of the random variables $(\boldsymbol{\theta}, \boldsymbol{\beta}, \boldsymbol{\varepsilon})$. Indeed, solving model $\widehat{\text{DEQ}}$ based on estimated weights $\hat{\underline{\omega}}_p$ that largely deviate from the ground truth $\underline{\omega}_p$ could lead to a severely suboptimal solution to the original CBCFLP and a poor approximation of its optimal value.

A natural measure for evaluating the discrepancy between the objective functions of models $\widehat{\text{DEQ}}$ and DEQ is the expected $L_1$ distance $\Phi(\boldsymbol{\omega}) := \mathbb{E}\left[\lvert\lvert\hat{\boldsymbol{\omega}}-\boldsymbol{\omega}\rvert\rvert_1\right]$. This measure corresponds to twice the expected total variation distance $\delta(\hat{W}, W)$ between the random variables $\hat{W}$ and $W$ with p.m.f.s $\hat{\boldsymbol{\omega}}$ and $\boldsymbol{\omega}$. We consider the normalized vector $\hat{\boldsymbol{\omega}}$ obtained through (\ref{eq:omega_hat}), with distribution $\hat{\boldsymbol{\omega}} \sim \text{Mult}(|R'|,\boldsymbol{\omega})$, where $|R'| = |R| - \sum_{r \in R}\mathds{1}[a_d(\boldsymbol{\theta}^r,\boldsymbol{\beta}^r,\boldsymbol{\varepsilon}^r)=0, \ \forall d \in D]$ is the number of non-trivial preferences profiles in the sample. The $L_1$ distance $\Phi(\boldsymbol{\omega})$ can be approximated as follows:
\begin{align*}
    \Phi(\boldsymbol{\omega}) &= \mathbb{E}\left[\lvert\lvert\hat{\boldsymbol{\omega}}-\boldsymbol{\omega}\rvert\rvert_1\right], \\
    &=\mathbb{E}\left[\sum_{p \in P} \left|\hat{\omega}_p-\omega_p\right|\right], \\
    &=\sum_{p \in P}  \mathbb{E}\left[\left|\hat{\omega}_p-\omega_p\right|\right],\\
    &\xrightarrow{|R'| \rightarrow \infty} \sum_{p \in P} \sqrt{\frac{2}{\pi}} \sqrt{\frac{\omega_p(1-\omega_p)}{|R'|}},\\
    &= \sqrt{\frac{2}{\pi |R'| }} \sum_{p \in P} \sqrt{\omega_p(1-\omega_p)},\\
    &=: \tilde{\Phi}(\boldsymbol{\omega}).
\end{align*}
\rev{The limit} is obtained by applying the approximation of the expected absolute error of the binomial parameter's estimator presented in \citet{blyth1980expected} to each term in the summation.

The following result allows us to draw a formal connection between the entropy of $W$ and the expected discrepancy between models $\widehat{\text{DEQ}}$ and DEQ. 

\begin{proposition}\label{prop:max_L1_error_general}
    For any parameter $\eta \in [0,1]$, the optimal solution of the following maximization problem:
    \begin{align}
        \label{inductiv_mod1}& \max_{\boldsymbol{\omega} \in [0, \eta]^{|P|}} \sum_{p \in P} \sqrt{\omega_p(1-\omega_p)}\\
        \label{inductiv_mod2} &\text{s.t. }   \sum_{p \in P} \omega_p = \eta,
    \end{align}
    is $\boldsymbol{\omega}^* = (\eta/|P|,\dots,\eta/|P|)$, with objective value $\sqrt{\eta(|P|-\eta)}$.
\end{proposition}
\proof
    The proof is by induction on $|P|$. The base case $|P|=1$ trivially holds, as $\boldsymbol{\omega}^* = (\eta)$ is the only feasible solution. Its objective value is $\sum_{p \in P}\sqrt{\omega^*_p(1-\omega^*_p)} = \sqrt{\eta(1-\eta)} = \sqrt{\eta(|P|-\eta)}$.

    Assuming that the result holds for a set of cardinality $|P| = k$, we demonstrate the case $|P| = k+1$. To do so, we consider an element $\bar{p} \in P$ and the set $\tilde{P} = P \setminus \{\bar{p}\}$, with $|\tilde{P}| = k$.
    \begin{align*}
        & \max_{\boldsymbol{\omega} \in [0, \eta]^{|P|}} \left\{\sum_{p \in P} \sqrt{\omega_p(1-\omega_p)} \Bigg| \sum_{p \in P} \omega_p = \eta \right\}\\
        =& \max_{\boldsymbol{\omega} \in [0, \eta]^{|P|}} \left\{ \sqrt{\omega_{\bar{p}}(1-\omega_{\bar{p}})} + \sum_{p \in \tilde{P}} \sqrt{\omega_p(1-\omega_p)} \Bigg|  \omega_{\bar{p}} + \sum_{p \in \tilde{P}} \omega_p = \eta \right\}\\
        =& \max_{\omega_{\bar{p}} \in [0, \eta]} \left\{ \sqrt{\omega_{\bar{p}}(1-\omega_{\bar{p}})} + \max_{\boldsymbol{\omega} \in [0, \eta - \omega_{\bar{p}}]^{k}} \left\{\sum_{p \in \tilde{P}} \sqrt{\omega_p(1-\omega_p)} \Bigg|\sum_{p \in \tilde{P}} \omega_p = \eta - \omega_{\bar{p}} \right\}\right\}\\
        =& \max_{\omega_{\bar{p}} \in [0, \eta]} \left\{ \sqrt{\omega_{\bar{p}}(1-\omega_{\bar{p}})} + \sqrt{(\eta - \omega_{\bar{p}})(k-(\eta - \omega_{\bar{p}}))} \right\}
    \end{align*}
    \rev{The last equality} is obtained by applying the inductive hypothesis to set $\tilde{P}$ with parameter $\eta-\omega_{\bar{p}}$. The inductive hypothesis also stipulates that the optimal solution to this inner problem is given by $\omega^*_{p}=(\eta-\omega_{\bar{p}})/|\tilde{P}| = (\eta-\omega_{\bar{p}})/k, \ \forall p \in \tilde{P}$. We now seek the maximizer $\omega^*_{\bar{p}}\in [0, \eta]$ of function $g(\omega_{\bar{p}}) := \left\{ \sqrt{\omega_{\bar{p}}(1-\omega_{\bar{p}})} + \sqrt{(\eta - \omega_{\bar{p}})(k-(\eta - \omega_{\bar{p}}))} \right\}$. 
    \begin{align}
    \label{eq:induction_critical1}\frac{\partial g(\omega_{\bar{p}})}{\partial \omega_{\bar{p}}} = 0  \iff& \frac{1-2\omega_{\bar{p}}}{2\sqrt{\omega_{\bar{p}}(1-\omega_{\bar{p}})}} + \frac{-k+2\eta-2\omega_{\bar{p}}}{2\sqrt{(\eta-\omega_{\bar{p}})(k-(\eta-\omega_{\bar{p}}))}} = 0\\
    \label{eq:induction_critical2}\iff& \omega_{\bar{p}} = \frac{\eta}{k+1}
    \end{align}

    Equation~(\ref{eq:induction_critical2}) is obtained through simple operations by developing expression (\ref{eq:induction_critical1}) and isolating $\omega_{\bar{p}}$. The only critical point $\omega^*_{\bar{p}} = \eta/(k+1)$ lies in the feasible interval $[0,\eta]$. The second derivative of function $g(\omega_{\bar{p}})$ is given by:

    \begin{align*}\frac{\partial^2 g(\omega_{\bar{p}})}{\partial \omega_{\bar{p}}^2} &= -\frac{k^2\sqrt{(\eta - \omega_{\bar{p}})(k-\eta+\omega_{\bar{p}})}}{4(\eta-\omega_{\bar{p}})^2(k-\eta+\omega_{\bar{p}})^2} - \frac{\sqrt{\omega_{\bar{p}}(1-\omega_{\bar{p}})}}{4\omega_{\bar{p}}^2(1-\omega_{\bar{p}})^2} <0.
    \end{align*}

    The maximizer of $g(\omega_{\bar{p}})$ is thus $\omega_{\bar{p}} = \eta/(k+1)$. It follows that $\omega^*_{p}=(\eta-\omega_{\bar{p}})/k = (\eta-\eta/(k+1))/k = \eta/(k+1) = \eta/|P|, \ \forall p \in \tilde{P}$. The optimal solution to problem (\ref{inductiv_mod1})-(\ref{inductiv_mod2}) in the inductive case is thus, as expected, $\boldsymbol{\omega}^* = (\eta/|P|,\dots,\eta/|P|)$.

    The optimal value, in accordance with the expected result, is given by:
    \begin{align*}
        \sum_{p \in P}\sqrt{\omega^*_p(1-\omega^*_p)} &= \sum_{p \in P}\sqrt{\frac{\eta}{|P|}(1-\frac{\eta}{|P|})},\\
        &= \sqrt{\eta(|P| - \eta)}.
    \end{align*}
\endproof

A direct consequence of Proposition~\ref{prop:max_L1_error_general} is that the maximum entropy distribution is also the maximizer of $\tilde{\Phi}(\boldsymbol{\omega})$. This is shown in the following corollary.

\begin{corollary}\label{corol:max_L1_error}
The estimated expected $L_1$ distance $\tilde{\Phi}(\boldsymbol{\omega})$ is maximized by the maximum entropy distribution $\boldsymbol{\omega}^* = (1/|P|, 1/|P|, \dots, 1/|P|)$.
\end{corollary}
\proof
We consider the following maximization problem.
\begin{align*}
    &\argmax_{\boldsymbol{\omega} \in [0, 1]^{|P|}} \left\{\tilde{\Phi}(\boldsymbol{\omega}) \Bigg| \sum_{p \in P} \omega_p = 1 \right\} \\
    =& \argmax_{\boldsymbol{\omega} \in [0, 1]^{|P|}} \left\{\sqrt{\frac{2}{\pi |R'| }} \sum_{p \in P} \sqrt{\omega_p(1-\omega_p)} \Bigg| \sum_{p \in P} \omega_p = 1\right\}\\
    \label{mod:max_L1_error}=& \argmax_{\boldsymbol{\omega} \in [0, 1]^{|P|}} \left\{\sum_{p \in P} \sqrt{\omega_p(1-\omega_p)} \Bigg| \sum_{p \in P} \omega_p = 1\right\}
\end{align*}
The optimal solution $\boldsymbol{\omega}^* = (1/|P|, 1/|P|, \dots, 1/|P|)$ to this last problem is given by Proposition~\ref{prop:max_L1_error_general}, with $\eta = 1$. The maximum estimated expected $L_1$ distance between $\boldsymbol{\hat{\omega}}$ and $\boldsymbol{\omega}$ is also given by Proposition~\ref{prop:max_L1_error_general}:
\begin{equation*}
    \tilde{\Phi}(\boldsymbol{\omega}^*) = \sqrt{\frac{2}{\pi |R'| }} \sqrt{1(|P|-1)} = \sqrt{\frac{2(|P|-1)}{\pi |R'| }}.
\end{equation*}
\endproof

Conversely, in the least entropy setting where $\omega_{p'} = 1 - \epsilon$ for a profile $p' \in P$, each term $\omega_p(1-\omega)$ converges to zero as $\epsilon \rightarrow 0$. \rev{For a fixed number of scenarios, we expect that the model obtained by sample average approximation is a better estimate of the distribution of the demand if a small set of profiles concentrates most of the demand in the original problem. In this section, we have shown that the measure of discrepancy $\tilde{\Phi}(\boldsymbol{\omega})$ between the original problem DEQ and its approximation $\widehat{\text{DEQ}}$ supports this conclusion.}

\subsection{Impact of Entropy on Computational Performance}\label{sec:info_theo:performance}
We expect \rev{PBD} to offer a significant advantage over SAAA \rev{for solving model $\widehat{\text{DEQ}}$} when the optimal value $\delta^*$ to the distance maximization problem (\ref{eq:knee}) is large. Indeed, this value corresponds to the difference between the proportion of the demand \rev{retained in the initial master problem} and the proportion of the decision variables $y_p$, $p \in \hat{P}$, included in the \rev{partial Benders reformulation}. A large value of $\delta^*$ thus means that \rev{the initial master problem solved in PBD is significantly smaller than model $\widehat{\text{DEQ}}$ but still takes most of the demand into account}. Hence, it limits the number of submodular cuts that have to be generated in the B\&C while significantly decreasing the number of decision variables compared to SAAA.

Once again, a formal connection between the entropy of the preference profiles and the computational properties of \rev{PBD} can be established. The following proposition serves this purpose. It shows that $\delta^*$ corresponds to the total variation distance between the empirical distribution of the observed preference profiles and the uniform discrete distribution $U$ defined on the same support $\hat{P}_1$.
\begin{proposition}
\label{prop:delta_tvd}
The maximum value $\delta^* = \max_{i \in \{1,\dots,|\hat{P}|\}} \Omega_i - i/|\hat{P}|$ is given by $\delta^* = \delta(\hat{W}, U) = 1/2 \lvert\lvert \hat{\boldsymbol{\omega}} - \bar{\boldsymbol{\omega}} \rvert\rvert_1$, where $\bar{\boldsymbol{\omega}} = \{1/|\hat{P}|,\dots,1/|\hat{P}|\}$ is the p.m.f. of the uniform discrete distribution $U$ over $\hat{P}$.
\end{proposition}
\proof
The difference between two subsequent values of $\delta_{i}$ is given by:
\begingroup
\allowdisplaybreaks
\begin{align*}
\delta_{i} - \delta_{i-1} &= \Omega_{i} - \frac{i}{|\hat{P}|} - \left( \Omega_{i-1} - \frac{i-1}{|\hat{P}|} \right),\\
&= \Omega_{i} - \Omega_{i-1} - \frac{1}{|\hat{P}|},\\
&= \sum_{k = 1}^{i} \hat{\omega}_{p_k} - \sum_{k = 1}^{i-1} \hat{\omega}_{p_k} - \frac{1}{|\hat{P}|},\\
&= \hat{\omega}_{p_i} - \frac{1}{|\hat{P}|}.
\end{align*}
\endgroup
Since $\{\hat{\omega}_{p_i}\}_{i=1}^{|\hat{P}|}$ is an increasing sequence, it follows from the definition of the knee index $i^*$ that $\hat{\omega}_{p_i} \geq \frac{1}{|\hat{P}|} \ \forall i \leq i^*$ and that $\hat{\omega}_{p_i} \leq \frac{1}{|\hat{P}|} \ \forall i > i^*$. The total variation distance $\delta(\hat{W}, U)$ can thus be developed as:
\begingroup
\allowdisplaybreaks
\begin{align*}
\delta(\hat{W}, U) &= \frac{1}{2}\lvert\lvert \hat{\boldsymbol{\omega}} - \bar{\boldsymbol{\omega}} \rvert\rvert_1,\\
&= \frac{1}{2}\sum_{k=1}^{|\hat{P}|} \left|\hat{\omega}_{p_k} - \frac{1}{|\hat{P}|}\right|,\\
&= \frac{1}{2} \left(\sum_{k=1}^{i^*} \left(\hat{\omega}_{p_k} - \frac{1}{|\hat{P}|}\right) + \sum_{k=i^*+1}^{|\hat{P}|} \left(\frac{1}{|\hat{P}|} - \hat{\omega}_{p_k}\right)\right),\\
&= \frac{1}{2} \left(\sum_{k=1}^{i^*} \hat{\omega}_{p_k} - \frac{i^*}{|\hat{P}|} + \frac{|\hat{P}| - i^*}{|\hat{P}|} - \sum_{k=i^*+1}^{|\hat{P}|}\hat{\omega}_{p_k}\right),\\
&= \frac{1}{2} \left(\Omega^*-\frac{i^*}{|\hat{P}|} + \frac{|\hat{P}|-i^*}{|\hat{P}|} - (1 - \Omega^*)\right),\\
&= \Omega^*-\frac{i^*}{|\hat{P}|},\\
&= \delta^*.
\end{align*}
\endgroup
\endproof

The total variation distance $\delta^* = \delta(\hat{W}, U)$ is linked to the entropy of $\hat{W}$ through the Kullback–Leibler divergence from $\hat{W}$ to $U$. From Pinsker's inequality \citep{pinsker1964information} and by Proposition~\ref{prop:delta_tvd}:
\begin{equation*}
    \delta^* \leq \sqrt{\frac{1}{2} D_{KL}(\hat{W}||U)} = \sqrt{\frac{1}{2} \left(\log|\hat{P}| - H(\hat{W})\right)}.
\end{equation*}

Other theoretical results, such as Bretagnolle–Huber inequality \citep{bretagnolle1978estimation},
\begin{equation*}
    \delta^* \leq \sqrt{1-e^{-D_{KL}(\hat{W}||U)}} = \sqrt{1-e^{H(\hat{W}) - \log|\hat{P}|}},
\end{equation*}
also provide a bound on the total variation distance based on the entropy of $\hat{W}$. The optimal value $\delta^*$ of the knee detection method used in the \rev{profile-retention strategy of PBD can be interpreted as the highest possible difference between the proportion of the total demand retained in the master problem and the proportion of profiles to which it corresponds. This value is thus bounded by monotonically decreasing functions of $H(\hat{W})$. The core idea of PBD is to take advantage of a small subset of profiles concentrating a high proportion of the demand in model $\widehat{\text{DEQ}}$. Consequently, for instances in which the entropy of $\hat{W}$ is close to its upper bound $\log|\hat{P}|$, PBD is less likely to offer a significant computational advantage over a general-purpose solver. An example illustrating the relative performance of \rev{PBD} and SAAA on instances of different levels of entropy is provided in Online Appendix~C. In this example, we also report computational results for the pure \rev{Benders decomposition} method obtained by setting $\hat{P}_1 = \emptyset$ in model $\widehat{\text{DEQ}}(\hat{P}_1)$ and give insights on why this approach generally performs poorly in the context of this work.} 

\subsection{\rev{Implications for Other Optimization Problems with Independent Agents}}\label{sec:info_theo:discussion}
\rev{Although our focus in this work is the CBCFLP and the instances of the MCLP that arise from simulation-based methods in this context, we want to highlight that the entropy measure and the results developed in this section apply to a broader class of problems. From a general perspective, a preference profile $p \in P$ can be defined as a function mapping each feasible strategic decision (here, a set of open facilities) to the reaction of an agent (here, the selection or not of a facility of the firm by a customer). The only assumption that was made throughout Section~\ref{sec:info_theo} is that the set $P$ is finite and that the weight $\omega_p$ of each profile $p \in P$ can be approximated by simulation, hence allowing to define a random variable $\hat{W}$ whose p.m.f. approximates that of the random variable $W$ with p.m.f. $\boldsymbol{\omega}$. As a consequence, our entropy measure can be used to quantify the level of concentration of the demand in any optimization problem in which independent agents reacting to strategic decisions can be partitioned into a finite set of profiles that completely characterize their behavior. In addition, it follows that the profile-retention strategy developed in Section~\ref{sec:simul_method:knee} could also be exploited in a partial Benders decomposition approach to such problems. Several deterministic and choice-based optimization problems have this structure.}

\rev{For example, in choice-based assortment optimization problems \citep[e.g.,][]{rusmevichientong2014assortment, zhang2020assortment}, the strategic decision is to select a subset of products, each associated with a fixed revenue, to offer to customers. Each customer then selects the available alternative with the highest utility, the opt-out alternative always being available. If the firm has access to $n$ products, the set of preference profiles $P$ thus corresponds to the $\left(n+1\right)!$ ways customers could rank the products and the opt-out by utility. The same information-theoretic analysis we presented in this section thus applies to such choice-based assortment optimization problems. If suitable Benders cuts were to be derived, so would a partial Benders decomposition method based on the Kneedle profile-retention strategy.}
 
\rev{On the contrary, our analysis would not directly apply to choice-based optimization problems involving continuous pricing of alternatives \citep[e.g.,][]{li2019product, Paneque_2021}. Indeed, in a very simple case involving a single pricing decision and a unique opt-out alternative, a preference profile would correspond to the maximum price at which a customer would be ready to select the alternative over the opt-out. If the utility function customers maximize comprises independent continuous random terms, as is the case under most commonly used RUMs, the set of preference profiles $P$ would have infinite cardinality. It follows that, for any finite number of draws, with probability one, each profile $p \in P$ would be observed at most if we were to approximate the distribution of $P$ by simulation, hence systematically leading $\hat{W}$ to follow the maximum entropy distribution on the support of the observed profiles $\hat{P}$.}

\rev{In contrast, for most deterministic covering problems, the demand can generally be represented through a small set of preference profiles. For example, \citet{cordeau2019benders} consider instances of the MCLP and of the partial set covering location problem in which a customer is covered by a candidate location if and only if their Euclidean distance is below a fixed threshold. In other words, each facility location covers a circular region around its position, and customers behave deterministically based on this rule. It is known that the maximal number of regions into which $n$ intersecting circles can divide the plane is $n^2-n+2$ \citep{gordon1987challenging}. All the instances considered in \citet{cordeau2019benders} include 100 candidate facility locations, meaning that at most 9,902 demand points, each corresponding to a region of the plane, would have to be considered after aggregating the customers by preference profile. The instances they study comprise a huge number of demand points, namely up to 15 and 40 million for the MCLP and the partial set covering location problem, respectively. Representing the demand based on preference profiles instead of individual customers would thus reduce the number of demand points of the largest instances of this testbed by at least four orders of magnitude. The entropy of the distribution of the resulting preference profiles would then be influenced by the location of the candidate facilities, their coverage radius, and the customers' spatial distribution. In a partial Benders decomposition approach to deterministic covering problems based on simple covering rules, the Kneedle profile-retention strategy would tend to retain in the master problem the profiles associated with the largest areas of the plane. These observations suggest that the theoretical analysis of this section and the partial Benders decomposition methodology developed in this work have relevant implications outside of the specific application considered in this paper.}

\section{Computational Experiments}\label{sec:experiments}
The purpose of this computational study is twofold. The first objective is to assess the potential of simulation-based methods for efficiently solving the CBCFLP under RUM models for which model-specific algorithms have already been extensively studied. To do so, \rev{PBD}, SAAA, and SAA are compared to GGX and MOA on MNL benchmark instances from the literature. The second objective is to compare the computational performance and solution quality provided by these methods under less restrictive modeling assumptions, which is done based on a new set of generative MMNL instances.

The experiments are conducted on a PC with processor Intel(R) i7-10875H CPU @ 2.30GHz along with 32 GB of RAM operated with Windows 10 Pro. The simulation-based methods \rev{PBD}, SAAA and SAA are implemented in Julia and are linked to IBM-ILOG CPLEX 20.1.0 optimization routines under default settings. The solver is warm-started with a simple greedy solution. The code and detailed computational results are available at \href{https://github.com/robinlegault/CBCFLP}{https://github.com/robinlegault/CBCFLP}. We use the original MATLAB implementation of MOA and GGX.

\subsection{Conditional MNL Instances}\label{sec:experiments:MNL}
The experiments of this section are carried out on two datasets on which the most recent methods for the CBCFLP under the MNL model have been benchmarked \citep{dam2022submodularity, Freire_2016, Ljubic_2018, Mai_2020}: 
\begin{itemize}
    \item HM14: This dataset has been proposed by \citet{Haase_2014}. To produce these problems, the location of the facilities and the customers have been uniformly generated on the square $[0,30]\times[0,30]$. This dataset includes instances with $|D|\in\{25,50,100\}$ available locations and $|N|\in\{50,100,200,400,800\}$ customers. In our experiments, we only consider the largest set of customers $|N|=800$.
    \item NYC: This dataset is based on a real-life park-and-ride location problem in New York City. It includes $|D|=59$ available locations and $|N|=\text{82,341}$ customers with weights ranging from 1 to 19. These are generally regarded as the most challenging MNL instances in the literature \citep{dam2022submodularity, Mai_2020}. \citet{holguin2012new} provide a more detailed presentation of the NYC instances.
\end{itemize}

For both datasets, the perceived utility of an open facility for customer $n \in N$ is given by $v^n_c = -\beta\theta^n_c$ for the new locations $c \in D$ and by $v^n_c = -\alpha\beta\theta^n_c$ for the competing locations $c \in E$. The attribute $\theta^n_{c}$ of customer $n \in N$ denotes its distance from facility $c \in C$ in HM14, and aggregates several factors (including travel time, tolls, auto costs, and waiting time) in NYC. The coefficient $\beta$ controls the importance of the deterministic term in the total utility and, thus, the entropy level of the preference profiles. The competitiveness of the existing facilities is controlled by $\alpha$. Smaller values of this coefficient lead to more attractive competitors. 

Like most previous studies, we consider a unique business constraint, namely that the firm can open at most $b$ facilities. Hence, the set of feasible configurations is given by $X = \{\boldsymbol{x} \in \{0,1\}^{|D|} : \sum_{d\in D} x_d \leq b\}$. As pointed out by \citet{Mai_2020} for the MNL model, opening additional locations cannot reduce the market share captured by the firm. This result trivially holds under any RUM model, as the effect of opening new facilities is to relax constraints (\ref{model:DEQ1:con}) in the deterministic equivalent model DEQ, or analogously constraints (\ref{model:DEQ_hat:con}) in model $\widehat{\text{DEQ}}$. An equality constraint can thus equivalently replace the inequality constraint in the definition of $X$.

For each number of available locations $|D| \in \{25, 50, 100\}$ in the HM14 dataset, we solve an instance for each configuration $(\beta, \alpha, b) \in \{1,2,5,10\} \times \{0.05, 0.1, 0.2\} \times \{2,3,4,5,6,7,8,9,10\}$. These parameters cover those used by \citet{Ljubic_2018} and \citet{Mai_2020}, except for $\alpha$. This parameter was taken in a wider interval in both studies, leading to degenerate instances in which more than 99.99\% of the market can be captured by a single location. The parameters we consider lead to more reasonable instances for which the optimal value ranges from 2\% to 90\% of the total market share. Similarly, each configuration $(\beta, \alpha, b) \in \{0.1, 0.15, 0.2, 0.25, 0.5, 0.75, 1, 1.25, 1.5, 1.75, 2\} \times \{0.75,1,1.25\} \times \{2,3,4,5,6,7,8,9,10\}$ is considered for the NYC dataset.

The simulation-based methods are applied with three different numbers of scenarios $|S_1|=10$, $|S_2|=100$ and $|S_3|=\text{1,000}$ for the HM14 instances and $|S_1|=1$, $|S_2|=5$ and $|S_3|=10$ for NYC. We use the code provided by the authors of the MOA algorithm, which sets $T=\min\{1000, |N|\}$ groups for the HM14 instances and $T=20$ for the NYC dataset.

The instances are grouped by number of candidate locations and value of $\beta$ to illustrate the impact of entropy on the performance of each method. \rev{For each group, we report the average empirical entropy $H(\hat{W})$ of the MCLP instances solved by the simulation-based methods and the upper bound $\log|\hat{P}|$ on this value for random variables whose support has the same cardinality. For a given family of instances, both the entropy and its upper bound tend to increase with $\beta$.} 

The average CPU times are reported in Table~\ref{table:MNL:time}. Table~\ref{table:MNL:model} presents the relative size of the models solved by SAAA and SAA (see columns $|\hat{P}|/|R|$) and by \rev{PBD} and SAAA (see columns $|\hat{P}_1|/|\hat{P}|$), as well as the optimal value $\delta^*$ of the knee detection problem. The last columns of Table~\ref{table:MNL:model} report the optimality gap for the simulation-based methods and GGX. For a feasible solution $x \in X$, it is defined as:
\begin{equation*}
    \text{Gap} = \frac{Z_N(\boldsymbol{x}^*_N) - Z_N(\boldsymbol{x})}{Z_N(\boldsymbol{x}^*_N)},
\end{equation*}
where $Z_N(\cdot)$ and $\boldsymbol{x}^*_N$ respectively denote the objective function and the optimal solution of the conditional MNL model~(\ref{eq:MNL_conditional}). The number of B\&C nodes explored by the simulation-based methods and the number of submodular cuts generated by \rev{PBD} are reported in Online Appendix~E.

\begin{table*}[htbp]
\centering
\caption{Average CPU times, in seconds, for conditional MNL instances, by entropy level (27 instances per row)}
\label{table:MNL:time}
\resizebox{\columnwidth}{!}{
\begin{tabular}{ccccrrrrrrrrrcr}
\toprule
\multirow{2}{*}{\makecell{Set}} & \multicolumn{3}{c}{Entropy} &  \multicolumn{3}{c}{\rev{PBD}} &  \multicolumn{3}{c}{SAAA} &  \multicolumn{3}{c}{SAA} & \multirowcell{2}{GGX} & \multirowcell{2}{MOA}\\
\cmidrule(lr){2-4}\cmidrule(lr){5-7}\cmidrule(lr){8-10}\cmidrule(lr){11-13} & $\beta$ & $H(\hat{W})$ & $\rev{\log|\hat{P}|}$ & $|S_1|$\pex  & $|S_2|$\pex & $|S_3|$\pex &$|S_1|$  & $|S_2|$ & $|S_3|$\pex & $|S_1|$  & $|S_2|$ & $|S_3|$\pex &  &   \\
\midrule 
\multirow{4}{*}{\makecell{HM14 \\ $|D|=25$}} 
   & 10 & 3.40 & \rev{\p 4.22} & 0.00\pex & 0.01\pex & 0.05\pex & 0.01 & 0.01 & 0.06\pex & 0.02 & 0.16 & 2.05\pex & 0.15 & 20.69 \\
   & 5  & 3.50 & \rev{\p 4.48} & 0.00\pex & 0.01\pex & 0.05\pex & 0.01 & 0.01 & 0.06\pex & 0.02 & 0.16 & 2.06\pex & 0.15 & 22.98 \\
   & 2  & 3.87 & \rev{\p 5.25} & 0.01\pex & 0.01\pex & 0.05\pex & 0.01 & 0.02 & 0.07\pex & 0.02 & 0.18 & 2.26\pex & 0.19 & 30.58 \\
   & 1  & 4.54 & \rev{\p 6.57} & 0.01\pex & 0.02\pex & 0.06\pex & 0.02 & 0.05 & 0.19\pex & 0.04 & 0.23 & 2.59\pex & 0.17 & 40.92 \\
\hline
\multirow{4}{*}{\makecell{HM14 \\ $|D|=50$}} 
   & 10 & 4.40 & \rev{\p 5.10} & 0.01\pex & 0.01\pex & 0.07\pex & 0.01 & 0.02 & 0.08\pex & 0.03 & 0.23 & 2.34\pex & 0.22 & 22.71 \\
   & 5  & 4.52 & \rev{\p 5.47} & 0.01\pex & 0.01\pex & 0.08\pex & 0.01 & 0.02 & 0.09\pex & 0.03 & 0.23 & 2.39\pex & 0.21 & 33.87 \\
   & 2  & 4.97 & \rev{\p 6.39} & 0.02\pex & 0.02\pex & 0.09\pex & 0.02 & 0.06 & 0.17\pex & 0.03 & 0.26 & 3.28\pex & 0.22 & 51.74 \\
   & 1  & 5.86 & \rev{\p 7.74} & 0.13\pex & 0.10\pex & 0.27\pex & 0.09 & 0.25 & 1.29\pex & 0.11 & 0.50 & 4.29\pex & 0.20 & 82.31 \\
\hline
\multirow{4}{*}{\makecell{HM14 \\ $|D|=100$}} 
    & 10 & 5.44 & \rev{\p 6.45} &  0.06\pex &  0.03\pex &  0.14\pex & 0.03 &  0.06 &  0.17\pex & 0.06 &  0.38 &  3.90\pex & 0.46 &  31.35 \\
    & 5  & 5.76 & \rev{\p 7.04} &  0.20\pex &  0.05\pex &  0.15\pex & 0.06 &  0.14 &  0.79\pex & 0.09 &  0.45 &  4.16\pex & 0.46 &  52.42 \\
    & 2  & 6.72 & \rev{\p 8.33} &  6.56\pex &  1.93\pex &  2.44\pex & 0.17 &  0.77 &  8.23\pex & 0.23 &  1.14 & 13.21\pex & 0.45 & 122.30 \\
    & 1  & 8.12 & \rev{\p 9.67} & 18.96$^{(8)}$ & 44.16$^{(8)}$ & 47.36$^{(9)}$ & 1.21 & 34.53 & 78.53$^{(8)}$ & 1.20 & 38.95 & 81.43$^{(8)}$ & 0.48 & 942.47 \\
\hline
\multirow{11}{*}{\makecell{NYC}} 
    & 2.00 & 3.88 & \rev{\p 6.17} & 0.02\pex &  0.06\pex &  0.10\pex & 0.07 &  0.27 &   0.44\pex & 0.30 &  1.63 &  3.33\pex & 5.95 & 22.44 \\
    & 1.75 & 3.91 & \rev{\p 6.26} & 0.02\pex &  0.06\pex &  0.10\pex & 0.07 &  0.29 &   0.46\pex & 0.30 &  1.57 &  3.33\pex & 6.01 & 17.12 \\
    & 1.50 & 3.96 & \rev{\p 6.35} & 0.02\pex &  0.06\pex &  0.10\pex & 0.08 &  0.31 &   0.49\pex & 0.30 &  1.59 &  3.43\pex & 5.93 & 11.95 \\
    & 1.25 & 4.04 & \rev{\p 6.51} & 0.02\pex &  0.06\pex &  0.11\pex & 0.09 &  0.37 &   0.61\pex & 0.31 &  1.65 &  3.53\pex & 5.95 &  7.25 \\
    & 1.00 & 4.17 & \rev{\p 6.76} & 0.02\pex &  0.07\pex &  0.12\pex & 0.10 &  0.43 &   0.77\pex & 0.32 &  1.80 &  3.79\pex & 5.94 &  4.49 \\
    & 0.75 & 4.42 & \rev{\p 7.14} & 0.02\pex &  0.09\pex &  0.15\pex & 0.19 &  0.66 &   1.14\pex & 0.42 &  2.10 &  4.48\pex & 5.93 &  2.70 \\
    & 0.50 & 4.99 & \rev{\p 7.89} & 0.03\pex &  0.12\pex &  0.21\pex & 0.32 &  1.32 &   2.59\pex & 0.50 &  2.94 &  6.49\pex & 5.95 &  1.54 \\
    & 0.25 & 6.62 & \rev{\p 9.47} & 0.06\pex &  0.32\pex &  0.68\pex & 1.11 &  6.09 &  12.00\pex & 1.42 &  8.85 & 20.42\pex & 5.95 &  0.75 \\
    & 0.20 & 7.25 & \rev{10.00} & 0.09\pex &  0.64\pex &  0.82\pex & 1.73 &  9.75 &  24.89\pex & 2.07 & 11.70 & 32.23\pex & 6.01 &  0.69 \\
    & 0.15 & 8.10 & \rev{10.66} & 0.14\pex &  1.02\pex &  2.79\pex & 3.44 & 19.95 &  41.21\pex & 3.32 & 21.12 & 57.15\pex & 5.92 &  0.60 \\
    & 0.10 & 9.37 & \rev{11.43} & 0.88\pex & 10.55\pex & 33.67\pex & 5.64 & 42.25 & 104.13\pex & 5.62 & 47.17 & 143.90\pex & 5.98 &  0.53 \\
\bottomrule
\multicolumn{14}{l}{\small{$^{(\cdot)}$ Number of instances that were not solved to optimality within the time limit of 10 minutes applied to simulation-based methods}}
\end{tabular}
}
\end{table*}
\begin{table*}[htbp]
\centering
\caption{Attributes of model $\widehat{\text{DEQ}}(\hat{P}_1)$ and solution quality compared to GGX for conditional MNL instances, by entropy level (27 instances per row)}
\label{table:MNL:model}
\resizebox{\columnwidth}{!}{
\begin{tabular}{ccccrrrrrrrrrcccc}
\toprule
\multirow{2}{*}{\makecell{Set}} & \multicolumn{3}{c}{Entropy} & \multicolumn{3}{c}{$|\hat{P}|/|R|$, (\%)} & \multicolumn{3}{c}{$|\hat{P}_1|/|\hat{P}|$, (\%)} & \multicolumn{3}{c}{$\delta^*$, (\%)} &  \multicolumn{4}{c}{{Gap, (\%)}} \\
\cmidrule(lr){2-4} \cmidrule(lr){5-7} \cmidrule(lr){8-10} \cmidrule(lr){11-13} \cmidrule(lr){14-17} & $\beta$ & $H(\hat{W})$ & $\rev{\log|\hat{P}|}$ &  $|S_1|$  & $|S_2|$ & $|S_3|$ &  $|S_1|$  & $|S_2|$ & $|S_3|$ & $|S_1|$  & $|S_2|$ & $|S_3|$ &  \p$|S_1|$  & $|S_2|$ & $|S_3|$\pexs  & GGX \\
\midrule 
\multirow{4}{*}{\makecell{HM14 \\ $|D|=25$}} 
   & 10 & 3.40 & \rev{\p 4.22} & 0.9 & 0.1 & 0.0 & 31.6 & 27.0 & 24.0 & 46.3 & 52.9 & 56.0 & \p.04 & .01 & \textbf{.00}\pexs & .10 \\
   & 5  & 3.50 & \rev{\p 4.48} & 1.1 & 0.2 & 0.0 & 27.2 & 21.1 & 18.5 & 49.4 & 56.4 & 60.4 & \p.12 & .02 & \textbf{.00}\pexs & \textbf{.00} \\
   & 2  & 3.87 & \rev{\p 5.25} & 2.0 & 0.4 & 0.1 & 22.0 & 16.9 & 12.2 & 52.1 & 64.9 & 74.5 & \p.18 & .03 & \textbf{.00}\pexs & \textbf{.00} \\
   & 1  & 4.54 & \rev{\p 6.57} & 4.7 & 1.3 & 0.3 & 19.0 & 11.5 &  7.8 & 61.9 & 76.5 & 85.2 & \p.31 & .11 & \textbf{.00}\pexs & \textbf{.00} \\
\hline
\multirow{4}{*}{\makecell{HM14 \\ $|D|=50$}} 
   & 10 & 4.40 & \rev{\p 5.10} & 2.1  & 0.3 & 0.0 & 32.9  & 31.2 & 28.5 & 40.9 & 47.0 & 52.5 &  \p.09 & \textbf{.00} & \textbf{.00}\pexs & \textbf{.00} \\
   & 5  & 4.52 & \rev{\p 5.47} & 2.7  & 0.4 & 0.1 & 30.7  & 25.4 & 20.8 & 45.0 & 55.1 & 62.7 &  \p.44 & \textbf{.00} & .01\pexs & .01 \\
   & 2  & 4.97 & \rev{\p 6.39} & 5.5  & 1.2 & 0.2 & 25.6  & 17.2 & 11.8 & 51.9 & 66.9 & 77.3 & 1.01 & .03 & \textbf{.00}\pexs & .02 \\
   & 1  & 5.86 & \rev{\p 7.74} & 12.9 & 4.3 & 1.3 & 20.6  & 13.2 &  8.6 & 55.6 & 72.5 & 83.1 & 1.11 & .08 & .01\pexs & \textbf{.00} \\
\hline
\multirow{4}{*}{\makecell{HM14 \\ $|D|=100$}} 
    & 10 & 5.44 & \rev{\p 6.45} &  6.8 &  1.2 &  0.2 & 29.2 & 23.9 & 18.7 & 45.6 & 57.0 & 64.9 &  \p.28 & .02 & \textbf{.00}\pexs & .40 \\
    & 5  & 5.76 & \rev{\p 7.04} & 10.4 &  2.3 &  0.5 & 24.3 & 18.7 & 13.7 & 48.5 & 62.5 & 72.0 &  \p.20 & .11 & \textbf{.00}\pexs & .16 \\
    & 2  & 6.72 & \rev{\p 8.33} & 22.0 &  8.6 &  2.9 & 18.2 & 15.4 & 10.3 & 46.6 & 65.6 & 78.5 & 1.11 & .15 & \textbf{.01}\pexs & .08 \\
    & 1  & 8.12 & \rev{\p 9.67} & 41.7 & 25.6 & 14.3 & 13.6 & 16.2 &  9.0 & 37.1 & 55.1 & 69.3 &  \p.62 & .11 & \textbf{.02}$^*$ & .06 \\
\hline
\multirow{11}{*}{\makecell{NYC}} 
& 2.00 & 3.88 & \rev{\p 6.17} &  3.2 &  1.3 &  0.9 & 12.9 & 10.2 &  9.4    & 70.5 & 76.3 & 78.0 & \p\textbf{.00} & \textbf{.00} & \textbf{.00}\pexs & \textbf{.00}  \\
& 1.75 & 3.91 & \rev{\p 6.26} &  3.3 &  1.4 &  1.0 & 12.3 & 10.1 &  9.0 & 70.6 & 76.7 & 78.4 & \p\textbf{.00} & \textbf{.00} & \textbf{.00}\pexs & \textbf{.00}  \\
& 1.50 & 3.96 & \rev{\p 6.35} &  3.6 &  1.6 &  1.1 & 12.5 &  9.7 &  8.9  & 70.7 & 76.9 & 78.3 & \p.01 & \textbf{.00} & \textbf{.00}\pexs & \textbf{.00}  \\
& 1.25 & 4.04 & \rev{\p 6.51} &  3.9 &  1.8 &  1.3 & 11.3 &  9.2 &  8.8 & 71.0 & 77.3 & 79.1 & \p\textbf{.00} & \textbf{.00} & \textbf{.00}\pexs & \textbf{.00}  \\
& 1.00 & 4.17 & \rev{\p 6.76} &  4.5 &  2.2 &  1.6 & 10.9 &  9.0 &  8.3    & 71.5 & 77.0 & 79.9 & \p\textbf{.00} & \textbf{.00} & \textbf{.00}\pexs & \textbf{.00}  \\
& 0.75 & 4.42 & \rev{\p 7.14} &  5.4 &  3.0 &  2.3 & 10.0 &  9.1 &  8.3 & 70.2 & 77.9 & 80.6 & \p\textbf{.00} & \textbf{.00} & \textbf{.00}\pexs & \textbf{.00}  \\
& 0.50 & 4.99 & \rev{\p 7.89} &  7.8 &  4.9 &  3.9 & 10.2 &  9.0 &  7.7  & 70.9 & 77.7 & 80.2 & \p.02 & .01 & .01\pexs & \textbf{.00}  \\
& 0.25 & 6.62 & \rev{\p 9.47} & 15.4 & 11.6 & 10.2 & 13.0 &  8.9 &  9.2 & 63.4 & 71.4 & 74.3 & \p.03 & \textbf{.00} & \textbf{.00}\pexs & \textbf{.00}  \\
& 0.20 & 7.25 & \rev{10.00} & 19.6 & 15.3 & 13.8 & 11.3 & 11.7 &  7.8  & 60.2 & 67.8 & 70.6 & \p\textbf{.00} & \textbf{.00} & \textbf{.00}\pexs & \textbf{.00}  \\
& 0.15 & 8.10 & \rev{10.66} & 26.7 & 21.9 & 20.1 &  9.8 &  9.8 & 10.1 & 54.7 & 61.8 & 64.6 & \p.01 & \textbf{.00} & \textbf{.00}\pexs & \textbf{.00}  \\
& 0.10 & 9.37 & \rev{11.43} & 40.9 & 35.4 & 33.3 & 16.4 & 17.3 & 17.5  & 45.5 & 51.5 & 53.8 & \p.01 & .03 & \textbf{.00}\pexs & \textbf{.00}  \\
\bottomrule
\multicolumn{16}{l}{\small{$^{*}$ Based on the instances that were solved to optimality within the time limit}}
\end{tabular}}
\end{table*}

The results in Table~\ref{table:MNL:time} indicate that \rev{PBD} is almost always the fastest among the simulation-based methods. An exception is the highest entropy instances of the HM14 dataset, for which all the simulation-based methods, especially \rev{PBD}, have relatively large computing times. SAAA consistently dominates SAA, most noticeably for instances with low entropy and a high number of simulated customers, for which the aggregation of simulated customers with identical preference profiles produces the most drastic reduction in the problem size, as reported in Table~\ref{table:MNL:model}. 

As expected, the computing time of \rev{PBD} is less sensitive to the number of scenarios than SAAA and SAA. In some cases, the CPU time of \rev{PBD} even decreases with the number of scenarios. This happens for the largest instances of the HM14 dataset with $\beta=2$, where \rev{PBD} requires 6.55 seconds on average for 10 scenarios and only 2.44 seconds for 1,000 scenarios. For the same group of instances, multiplying the number of scenarios by 100 (comparing $|S_1|$ and $|S_3|$) increases the CPU time of SAAA and SAA by a factor of 48 and 58, respectively. Supporting the theoretical analysis of Section~\ref{sec:info_theo:performance}, the progressive improvement in the relative performance of \rev{PBD} compared to SAAA and SAA that comes with more scenarios is accompanied by an increase in the optimal value $\delta^*$ of the knee detection problem (see Table~\ref{table:MNL:model}). In the previously discussed group with 10 scenarios, $\hat{P}_1$ comprises 18.2\% of the observed preference profiles accounting for 65.1\% of the demand, for a difference of $\delta^*=46.6\%$. In comparison, with 1,000 scenarios, $\delta^*=78.5\%$ and the 10.3\% most important profiles account for 88.8\% of the demand. Furthermore, across both datasets, \rev{PBD} requires exploring fewer nodes and generating fewer submodular cuts in the B\&C tree as the number of scenarios increases. SAAA and SAA solve almost all the instances to optimality at the root node (see Online Appendix~E). 

Since GGX only relies on local search procedures, its execution time is mostly unaffected by the entropy level, whereas the instances with high entropy are generally more challenging for the other methods. The only exception occurs in the NYC dataset, where the CPU time of MOA decreases significantly with the entropy. A possible explanation is the large number of customers (4117 on average for $T=20$) aggregated in each component $g_t$ of the objective function. When the entropy is low, the customer's choice probabilities are mainly determined by their attributes instead of the random error term. As $\beta$ increases, so does the heterogeneity of the customers' preferences in each group, which we conjecture makes the aggregated subgradient cuts weaker and negatively affects the performance of this version of MOA. This behavior is not observed on the HM14 dataset, as this set of instances is solved using the multicut version of MOA.

The results of Table~\ref{table:MNL:model} support the analysis of Section~\ref{sec:info_theo:quality}, as they illustrate the negative impact of entropy on the solution quality obtained by sample average approximation. This is especially visible for the smallest number of scenarios on the HM14 dataset, where the optimality gap is, on average, five times larger with $\beta=1$ than with $\beta=10$. Fortunately, closing the optimality gap does not require an unmanageable number of scenarios, even for high entropy instances. For the largest number of scenarios, \rev{PBD} remains significantly faster than MOA and even GGX in most cases and provides near-optimal solutions for all the instances, with an average optimality gap that does not exceed 0.02\%. In comparison, the optimality gap of GGX reaches 0.40\% for the HM14 dataset with $|D|=100$ and $\beta=10$. For the same group of instances, \rev{PBD} with 1,000 scenarios requires 229 times less CPU time than MOA and three times less than GGX and leads to a negligible optimality gap. 

These results demonstrate that \rev{PBD}, despite being a model-free approach, can provide better solution quality and computational performance than the state-of-the-art heuristic method for conditional MNL instances. In addition, its asymptotic optimality property makes it a sound alternative to exact methods for challenging large-scale instances.

\subsection{Generative MMNL Instances} 
\label{sec:experiments:MMNL}
This section presents experiments based on a new set of generative MMNL instances. In the MIX dataset, we consider three types of locations and three types of customers with different utility functions. The customers are distributed across four neighborhoods with different population densities and demographic properties. These instances include $|E|=30$ competing facilities and $|D|=75$ available locations. A random multiplier $\beta \sim U[\beta^-, \beta^+]$ controls the weight of the perceived utility $ v_c(\boldsymbol{\theta})$ in the utility function $u_c(\boldsymbol{\theta}, \boldsymbol{\varepsilon}) = v_c(\boldsymbol{\theta}) + \varepsilon_c$. The level of entropy of the preference profiles in the population decreases with $\beta$. A detailed presentation of the MIX dataset is provided in Online Appendix~D. 

To solve the generative MMNL model with MOA and GGX, we sample a set $N$ of realizations $\{(\boldsymbol{\theta}^n, \beta^n)\}_{n \in N}$ of the customers' random attributes $\boldsymbol{\theta}$ and of the random coefficient $\beta$. For \rev{PBD}, SAAA and SAA, the resulting conditional MNL model is approximated as in the previous section through a set $S$ of realizations $\{\boldsymbol{\varepsilon}^{ns}\}_{s \in S}$ of the error component for each customer $n \in N$. 

For the same number of realizations, solving the MNL model exactly with MOA provides the best solution on average but is more expensive than solving it heuristically with GGX, which in turn is more expensive than solving the simulation-based model with \rev{PBD}. To account for the different computational performance and solution quality of each approach, different numbers of customers are considered for each of them. We generate $|N| \in \{125, 250, 500, 1000, 2000\}$ customers for MOA,  $|N| \in \{250, 1000, 4000, 16000, 64000\}$ for GGX, and $|N| \in \{16000, 32000, 64000, 128000, 256000\}$ with $|S|=5$ scenarios for \rev{PBD}, SAAA and SAA. We study the model with four different pairs of bounds $(\beta^-, \beta^+) \in \{(0.125,0.25), (0.25,0.5), (0.5,1), (1,2)\}$, which control the level of stochasticity of the RUM model. For each set of parameters and each sample size, five instances are generated and solved with five different budgets $b \in \{5, 10, 15, 20, 25\}$. \rev{We report the average empirical entropy $H(\hat{W})$ of the groups of instances solved by the simulation-based methods and the associated upper bound $\log|\hat{P}|$.}

Table~\ref{table:MMNL:time} reports the CPU times for each method and the relative size of the different simulation-based models. For a sample of customers $N$ and a set of scenarios $S$, let us respectively denote by $\boldsymbol{x}^*_{NS}$, $\boldsymbol{x}^{G}_{N}$ and $\boldsymbol{x}^*_{N}$ the optimal solution of the simulation-based model $\widehat{\text{DEQ}}$, the heuristic solution returned by GGX, and the optimal solution of the conditional MNL model. Also, let $Z_{NS}(\cdot)$, $Z_{N}(\cdot)$ and $Z(\cdot)$ be the objective functions of the simulation-based, conditional MNL, and generative MMNL models. Table~\ref{table:MMNL:value} presents the average values of solutions $\boldsymbol{x}^*_{NS}$, $\boldsymbol{x}^G_{N}$ and $\boldsymbol{x}^*_{N}$ for each objective function. $Z(\cdot)$ is evaluated based on an independent sample of $|N|=\text{1,000,000}$ customers.

\begin{table*}[htbp]
\centering
\caption{Average CPU times and attributes of model $\widehat{\text{DEQ}}(\hat{P}_1)$ for the MIX dataset, by entropy level and sample size (25 instances per row)}
\label{table:MMNL:time}
\resizebox{\columnwidth}{!}{
\begin{tabular}{ccccccccccrrrrr}
\toprule
\multicolumn{4}{c}{Entropy} & \multicolumn{3}{c}{$|N|$, (thousands)} & \multicolumn{3}{c}{Model $\widehat{\text{DEQ}}(\hat{P}_1)$, (\%)} &  \multicolumn{5}{c}{Time, (seconds)} \\
\cmidrule(lr){1-4}\cmidrule(lr){5-7}\cmidrule(lr){8-10}\cmidrule(lr){11-15} $\beta^-$ & $\beta^+$ & $H(\hat{W})$ & $\rev{\log|\hat{P}|}$  & \rev{PBD} & GGX & MOA & $|\hat{P}|/|R|$ & $|\hat{P}_1|/|\hat{P}|$ & $\delta^*$ & \rev{PBD} & SAAA & SAA & GGX & MOA  \\
\midrule 
\multirow{5}{*}{1} & \multirow{5}{*}{2} 
   & 5.09 & \rev{\p 7.14} &  16 &  0.25 & 0.125  & \p1.6 & 13.6 & 74.2 &   0.02 &   0.03 &   0.33 &  0.73 &   0.43 \\
 & & 5.11 & \rev{\p 7.36} &  32 &  1    & 0.25   & \p1.0 & 12.6 & 77.1 &   0.03 &   0.05 &   0.65 &  1.20 &   1.27 \\
 & & 5.12 & \rev{\p 7.56} &  64 &  4    & 0.5    & \p0.6 & 11.5 & 79.4 &   0.05 &   0.07 &   1.76 &  7.40 &   7.50 \\
 & & 5.12 & \rev{\p 7.77} & 128 & 16    & 1      & \p0.4 & 10.3 & 81.7 &   0.10 &   0.12 &   3.81 & 24.26 &  56.59 \\
 & & 5.12 & \rev{\p 7.95} & 256 & 64    & 2      & \p0.2 &  \p9.6 & 83.4 &   0.18 &   0.19 &   7.27 & 47.17 & 126.83 \\
\hline\multirow{5}{*}{0.5} & \multirow{5}{*}{1} 
   & 5.61 & \rev{\p 7.91} &  16 &  0.25 & 0.125  &  \p3.4 & 11.2 & 75.7 & 0.03 &   0.06 &   0.36 &  0.73 &   0.47 \\
 & & 5.63 & \rev{\p 8.23} &  32 &  1    & 0.25   &  \p2.4 &  \p9.7 & 79.0 & 0.04 &   0.08 &   0.73 &  1.20 &   1.78 \\
 & & 5.65 & \rev{\p 8.54} &  64 &  4    & 0.5    &  \p1.6 &  \p8.6 & 81.7 & 0.06 &   0.12 &   1.75 &  7.69 &   9.82 \\
 & & 5.65 & \rev{\p 8.83} & 128 & 16    & 1      &  \p1.1 &  \p7.6 & 84.1 & 0.10 &   0.19 &   3.71 & 24.44 &  69.35 \\
 & & 5.67 & \rev{\p 9.11} & 256 & 64    & 2      &  \p0.7 &  \p6.9 & 86.0 & 0.20 &   0.30 &   8.06 & 47.29 & 135.05 \\
\hline\multirow{5}{*}{0.25} & \multirow{5}{*}{0.5} 
   & 6.63 & \rev{\p 8.93} &  16 &  0.25 & 0.125 &  \p9.5 & 10.0 & 71.2 &   0.07 &   0.28 &   0.63 &  0.71 &   1.02 \\
 & & 6.68 & \rev{\p 9.36} &  32 &  1    & 0.25  &  \p7.3 &  \p8.1 & 74.8 &   0.12 &   0.61 &   1.39 &  1.25 &   3.06 \\
 & & 6.74 & \rev{\p 9.79} &  64 &  4    & 0.5   &  \p5.6 &  \p7.4 & 77.9 &   0.16 &   1.02 &   3.07 &  7.45 &  14.76 \\
 & & 6.77 & \rev{  10.19} & 128 & 16    & 1     &  \p4.2 &  \p6.8 & 80.8 &   0.21 &   1.23 &   7.42 & 24.16 &  95.18 \\
 & & 6.80 & \rev{  10.59} & 256 & 64    & 2     &  \p3.1 &  \p6.1 & 83.2 &   0.34 &   7.15 &  18.59 & 47.25 & 182.12 \\
\hline\multirow{5}{*}{0.125} & \multirow{5}{*}{0.25} 
   & 7.86 & \rev{\p 9.77} &  16 &  0.25 & 0.125 & 21.9 &  \p8.7 & 59.5 &   0.49 &   1.64 &   1.87 &  0.73 &   1.74 \\
 & & 8.01 & \rev{  10.30} &  32 &  1    & 0.25  & 18.5 &  \p9.2 & 63.7 &   0.66 &   5.68 &   7.28 &  1.26 &   6.06 \\
 & & 8.16 & \rev{  10.83} &  64 &  4    & 0.5   & 15.9 &  \p7.7 & 67.4 &   0.78 &   6.87 &  11.66 &  7.17 &  20.70 \\
 & & 8.29 & \rev{  11.36} & 128 & 16    & 1     & 13.4 &  \p6.7 & 70.7 &   1.31 &  52.98 &  59.88 & 23.75 & 115.48 \\
 & & 8.39 & \rev{  11.89} & 256 & 64    & 2     & 11.3 &  \p6.2 & 73.7 &   2.44 & 134.20 & 218.34 & 48.03 & 206.64 \\
\bottomrule
\end{tabular}
}
\end{table*}

The results presented in Table~\ref{table:MMNL:time} indicate that \rev{PBD} is consistently the most efficient simulation-based method for the MIX dataset. In the lowest entropy settings, the advantage of \rev{PBD} over SAAA is only marginal, as the number of observed preference profiles remains limited even for a very large number of simulated customers. For example, for $(\beta^-, \beta^+)=(1,2)$ with $|N|=\text{256,000}$ customers and $|S|=5$ scenarios, aggregating the simulated customers by preference profile reduces the model's size from $|R|=\text{1,280,000}$ to $|\hat{P}|=\text{2,840}$ (i.e., $|\hat{P}|/|R| \approx 0.2\%$). However, for parameters $(\beta^-, \beta^+) = (0.125, 0.25)$ and the same number of simulated customers, the average number of observed preference profiles reaches $|\hat{P}|=\text{145,175}$ (i.e., $|\hat{P}|/|R| \approx 11.3\%$). This makes solving model $\widehat{\text{DEQ}}$ with CPLEX computationally demanding. In this setting, $\hat{P}$ contains a very high number of preference profiles with negligible weight, and aggregating their contribution to the objective value becomes increasingly profitable computationally. The subset $\hat{P}_1$ of preference profiles that are explicitly represented in model $\widehat{\text{DEQ}}(\hat{P}_1)$ has an average cardinality of 8,998 (6.2\% of the observed preference profiles) and accounts for approximately 80\% of the demand ($\Omega^* = |\hat{P}_1|/|\hat{P}|+\delta^* \approx 79.9\%$). Although $\hat{P}_2$ regroups 136,177 preference profiles on average, \rev{PBD} terminates after generating only 23.6 submodular cuts on average (see Online Appendix~E). As a result, the average CPU time of \rev{PBD} for this group of instances is 2.44 seconds compared to 134.20 seconds and 218.34 seconds for SAAA and SAA.

\begin{table*}[htbp]
\centering
\caption{Average objective value evaluated on an independent sample of 1,000,000 customers (expected market share, in \%) under the simulation-based, conditional MNL and generative MMNL models for the MIX dataset, by entropy level and sample size (25 instances per row)}
\label{table:MMNL:value}
\resizebox{\columnwidth}{!}{
\begin{tabular}{cccccccccccccc}
\toprule
\multicolumn{4}{c}{Entropy} &  \multicolumn{3}{c}{$|N|$, (thousands)} & \multicolumn{3}{c}{\rev{PBD}} & \multicolumn{2}{c}{GGX} & \multicolumn{2}{c}{MOA} \\
\cmidrule(lr){1-4}\cmidrule(lr){5-7}\cmidrule(lr){8-10}\cmidrule(lr){11-12}\cmidrule(lr){13-14} $\beta^-$ & $\beta^+$ & $H(\hat{W})$ & $\rev{\log|\hat{P}|}$ & \rev{PBD} & GGX & MOA & $Z_{NS}(\boldsymbol{x}^*_{NS})$ & $Z_{N}(\boldsymbol{x}^*_{NS})$ & $Z(\boldsymbol{x}^*_{NS})$ & $Z_{N}(\boldsymbol{x}^G_{N})$ & $Z(\boldsymbol{x}^G_{N})$ & $Z_{N}(\boldsymbol{x}^*_{N})$ & $Z(\boldsymbol{x}^*_{N})$ \\
\midrule 
\multirow{5}{*}{1} & \multirow{5}{*}{2} 
   & 5.09 & \rev{\p 7.14} &  16 &  0.25 & 0.125 & 62.59 & 62.54 & 62.62 & 63.38 & 61.02 & 65.59 & 60.29  \\
 & & 5.11 & \rev{\p 7.36} &  32 &  1    & 0.25  & 62.71 & 62.70 & 62.62 & 62.84 & 62.35 & 63.38 & 61.02  \\
 & & 5.12 & \rev{\p 7.56} &  64 &  4    & 0.5   & 62.49 & 62.50 & \textbf{62.63} & 62.78 & 62.51 & 63.55 & 61.32  \\
 & & 5.12 & \rev{\p 7.77} & 128 & 16    & 1     & 62.52 & 62.53 & 62.62 & 62.55 & 62.62 & 62.84 & 62.35  \\
 & & 5.12 & \rev{\p 7.95} & 256 & 64    & 2     & 62.59 & 62.59 & \textbf{62.63} & 62.50 & \textbf{62.63} & 63.02 & 62.44  \\
\hline\multirow{5}{*}{0.5} & \multirow{5}{*}{1} 
   & 5.61 & \rev{\p 7.91} &  16 &  0.25 & 0.125 & 60.63 & 60.58 & 60.70 & 60.98 & 59.28 & 63.68 & 58.71  \\
 & & 5.63 & \rev{\p 8.23} &  32 &  1    & 0.25  & 60.64 & 60.67 & 60.70 & 61.09 & 60.44 & 60.98 & 59.28  \\
 & & 5.65 & \rev{\p 8.54} &  64 &  4    & 0.5   & 60.65 & 60.68 & 60.70 & 60.81 & 60.68 & 61.25 & 60.16  \\
 & & 5.65 & \rev{\p 8.83} & 128 & 16    & 1     & 60.66 & 60.67 & \textbf{60.71} & 60.58 & 60.70 & 61.09 & 60.44  \\
 & & 5.67 & \rev{\p 9.11} & 256 & 64    & 2     & 60.72 & 60.74 & \textbf{60.71} & 60.68 & \textbf{60.71} & 61.35 & 60.61  \\
\hline\multirow{5}{*}{0.25} & \multirow{5}{*}{0.5} 
   & 6.63 & \rev{\p 8.93} &  16 &  0.25 & 0.125 & 56.96 & 56.88 & 56.87 & 58.37 & 55.90 & 58.69 & 55.04  \\
 & & 6.68 & \rev{\p 9.36} &  32 &  1    & 0.25  & 56.90 & 56.95 & 56.87 & 56.41 & 56.70 & 58.37 & 55.90  \\
 & & 6.74 & \rev{\p 9.79} &  64 &  4    & 0.5   & 56.86 & 56.84 & \textbf{56.88} & 56.87 & 56.86 & 57.17 & 56.33  \\
 & & 6.77 & \rev{  10.19} & 128 & 16    & 1     & 56.89 & 56.92 & \textbf{56.88} & 56.88 & 56.87 & 56.41 & 56.70  \\
 & & 6.80 & \rev{  10.59} & 256 & 64    & 2     & 56.92 & 56.92 & \textbf{56.88} & 56.84 & \textbf{56.88} & 57.52 & 56.81  \\
\hline\multirow{5}{*}{0.125} & \multirow{5}{*}{0.25} 
   & 7.86 & \rev{\p 9.77} &  16 &  0.25 & 0.125 & 51.43 & 51.29 & 51.23 & 53.05 & 50.57 & 52.24 & 49.89  \\
 & & 8.01 & \rev{  10.30} &  32 &  1    & 0.25  & 51.04 & 51.07 & 51.24 & 50.85 & 51.09 & 53.05 & 50.57  \\
 & & 8.16 & \rev{  10.83} &  64 &  4    & 0.5   & 51.24 & 51.27 & \textbf{51.25} & 51.28 & 51.21 & 52.23 & 50.94  \\
 & & 8.29 & \rev{  11.36} & 128 & 16    & 1     & 51.18 & 51.20 & \textbf{51.25} & 51.31 & \textbf{51.25} & 50.85 & 51.09  \\
 & & 8.39 & \rev{  11.89} & 256 & 64    & 2     & 51.25 & 51.24 & \textbf{51.25} & 51.28 & \textbf{51.25} & 51.16 & 51.18  \\
\bottomrule
\end{tabular}
}
\end{table*}

The objective value of solution $\boldsymbol{x}^*_{NS}$ for the generative MMNL model stabilizes at a near-optimal level between $|N|=\text{64,000}$ and $|N|=\text{256,000}$ (see figures highlighted in bold in Table~\ref{table:MMNL:value}). \rev{PBD} generally solves these instances in less than one second. Similar quality solutions can be obtained in approximately 47 seconds using GGX with $|N|=\text{64,000}$. MOA is the worst-performing method for the MIX dataset. Indeed, solving the conditional MNL model to optimality for $|N|=\text{2,000}$ customers is three to four orders of magnitude longer than solving the simulation-based model with \rev{PBD} for $|N|=\text{16,000}$ and $|S|=5$ and consistently provides solutions of lower quality for the generative MMNL model.

The results of this section demonstrate the clear advantage of our simulation-based method over model-specific algorithms and classical sample average approximation for problems based on flexible RUM models and a generative perspective.

\section{Conclusion}\label{sec:conclusion}
This paper presents a model-free approach for solving probabilistic competitive facility location problems with utility-maximizing customers. We introduce a deterministic equivalent reformulation of the problem, \rev{from which it follows that the objective function of the competitive facility location problem is submodular under any RUM model. Approximating this deterministic reformulation by simulation leads to a maximum covering location problem in which each demand point represents a different preference profile. We propose a partial Benders reformulation in which the preference profiles are partitioned into two sets that are respectively explicitly represented in the model and aggregated into a unique composite customer whose contribution to the objective is bounded by submodular cuts.}

\rev{Our partial Benders decomposition exploits a new profile-retention strategy that is based on a knee detection method and aims at including in the master problem a small set of preference profiles that account for a large part of the total demand.} As a result, our algorithm does not rely on user-defined parameters. We develop an information-theoretic analysis of the problem and draw connections between the entropy of the preference profiles and the computational performance of our approach. Computational experiments on conditional MNL and generative MMNL instances show that our method can significantly outperform state-of-the-art model-specific algorithms in terms of computing time and solution quality. A key takeaway from the experiments is that our methodology scales significantly better than the classical sample average approximation method with respect to the number of simulated customers. Combined with the absence of restrictive modeling assumptions of our approach, this opens the way to integrating larger and more complex populations in choice-based competitive facility location problems. 

Regarding future research directions, \rev{our partial Benders decomposition} method could be generalized to a multicut version including multiple auxiliary variables bounded by independent sets of cuts. This approach \rev{could benefit from a different profile-retention strategy} than the knee detection method used in our single-cut version, \rev{and would require a partitioning strategy for the profiles that are not retained in the master problem}. \rev{It would also be interesting to extend and apply our methodology to other choice-based optimization problems that can be approached by simulation, such as assortment optimization.} Relevant modeling extensions of this work include adding capacity constraints on the facilities and accounting for the reaction of the competitors to the firm's decisions.

\section*{Acknowledgments}
\label{sec:acknowledgments}
This research was supported by FRQNT and IVADO through scholarships to the first author. The second author \rev{was partially supported by the Canada Research Chair program [950-232244]}. We wish to thank Tien Mai and Andrea Lodi for sharing the code of the MOA and GGX algorithms as well as their benchmark instances. \rev{We are grateful for the two excellent reviewer reports that helped us greatly improve the quality of the manuscript and the positioning of our work.}

\bibliography{main}

\begin{thebibliography}{53}
\providecommand{\natexlab}[1]{#1}
\providecommand{\url}[1]{\texttt{#1}}
\expandafter\ifx\csname urlstyle\endcsname\relax
  \providecommand{\doi}[1]{doi: #1}\else
  \providecommand{\doi}{doi: \begingroup \urlstyle{rm}\Url}\fi

\bibitem[Aros-Vera et~al.(2013)Aros-Vera, Marianov, and Mitchell]{aros_2013}
Aros-Vera, F., Marianov, V., and Mitchell, J.~E.
\newblock p-{H}ub approach for the optimal park-and-ride facility location problem.
\newblock \emph{European Journal of Operational Research}, 226\penalty0 (2):\penalty0 277--285, 2013.

\bibitem[Benati(1997)]{benati1997submodularity}
Benati, S.
\newblock Submodularity in competitive location problems.
\newblock \emph{Ricerca Operativa}, 26:\penalty0 3--34, 1997.

\bibitem[Benati and Hansen(2002)]{Benati_2002}
Benati, S. and Hansen, P.
\newblock The maximum capture problem with random utilities: Problem formulation and algorithms.
\newblock \emph{European Journal of Operational Research}, 143\penalty0 (3):\penalty0 518--530, 2002.

\bibitem[Berbeglia and Joret(2020)]{berbeglia2020assortment}
Berbeglia, G. and Joret, G.
\newblock Assortment optimisation under a general discrete choice model: A tight analysis of revenue-ordered assortments.
\newblock \emph{Algorithmica}, 82\penalty0 (4):\penalty0 681--720, 2020.

\bibitem[Bhat and Guo(2004)]{Bhat_2004}
Bhat, C.~R. and Guo, J.
\newblock A mixed spatially correlated logit model: formulation and application to residential choice modeling.
\newblock \emph{Transportation Research Part B: Methodological}, 38\penalty0 (2):\penalty0 147--168, 2004.

\bibitem[Birge and Louveaux(1988)]{Birge_1988}
Birge, J.~R. and Louveaux, F.~V.
\newblock A multicut algorithm for two-stage stochastic linear programs.
\newblock \emph{European Journal of Operational Research}, 34\penalty0 (3):\penalty0 384--392, 1988.

\bibitem[Blyth(1980)]{blyth1980expected}
Blyth, C.~R.
\newblock Expected absolute error of the usual estimator of the binomial parameter.
\newblock \emph{The American Statistician}, 34\penalty0 (3):\penalty0 155--157, 1980.

\bibitem[Bonami et~al.(2008)Bonami, Biegler, Conn, Cornuéjols, Grossmann, Laird, Lee, Lodi, Margot, Sawaya, and Wächter]{Bonami_2008}
Bonami, P., Biegler, L.~T., Conn, A.~R., Cornuéjols, G., Grossmann, I.~E., Laird, C.~D., Lee, J., Lodi, A., Margot, F., Sawaya, N., and Wächter, A.
\newblock An algorithmic framework for convex mixed integer nonlinear programs.
\newblock \emph{Discrete Optimization}, 5\penalty0 (2):\penalty0 186--204, 2008.

\bibitem[Bretagnolle and Huber(1978)]{bretagnolle1978estimation}
Bretagnolle, J. and Huber, C.
\newblock Estimation des densit{\'e}s: risque minimax.
\newblock \emph{S{\'e}minaire de probabilit{\'e}s de Strasbourg}, 12:\penalty0 342--363, 1978.

\bibitem[Chu(1989)]{Chu_1989}
Chu, C.
\newblock A paired combinatorial logit model for travel demand analysis.
\newblock In \emph{Proceedings of the Fifth World Conference on Transportation Research, 1989}, volume~4, pages 295--309, 1989.

\bibitem[Church and ReVelle(1974)]{church1974maximal}
Church, R. and ReVelle, C.
\newblock The maximal covering location problem.
\newblock In \emph{Papers of the regional science association}, volume~32, pages 101--118. Springer-Verlag Berlin/Heidelberg, 1974.

\bibitem[Coniglio et~al.(2022)Coniglio, Furini, and Ljubi{\'c}]{coniglio2022submodular}
Coniglio, S., Furini, F., and Ljubi{\'c}, I.
\newblock Submodular maximization of concave utility functions composed with a set-union operator with applications to maximal covering location problems.
\newblock \emph{Mathematical Programming}, 196\penalty0 (1-2):\penalty0 9--56, 2022.

\bibitem[Cordeau et~al.(2019)Cordeau, Furini, and Ljubi{\'c}]{cordeau2019benders}
Cordeau, J.-F., Furini, F., and Ljubi{\'c}, I.
\newblock Benders decomposition for very large scale partial set covering and maximal covering location problems.
\newblock \emph{European Journal of Operational Research}, 275\penalty0 (3):\penalty0 882--896, 2019.

\bibitem[Crainic et~al.(2021)Crainic, Hewitt, Maggioni, and Rei]{crainic2021partial}
Crainic, T.~G., Hewitt, M., Maggioni, F., and Rei, W.
\newblock Partial benders decomposition: general methodology and application to stochastic network design.
\newblock \emph{Transportation Science}, 55\penalty0 (2):\penalty0 414--435, 2021.

\bibitem[Dam et~al.(2022)Dam, Ta, and Mai]{dam2022submodularity}
Dam, T.~T., Ta, T.~A., and Mai, T.
\newblock Submodularity and local search approaches for maximum capture problems under generalized extreme value models.
\newblock \emph{European Journal of Operational Research}, 300\penalty0 (3):\penalty0 953--965, 2022.

\bibitem[Davis et~al.(2017)Davis, Topaloglu, and Williamson]{davis2017pricing}
Davis, J.~M., Topaloglu, H., and Williamson, D.~P.
\newblock Pricing problems under the nested logit model with a quality consistency constraint.
\newblock \emph{INFORMS Journal on Computing}, 29\penalty0 (1):\penalty0 54--76, 2017.

\bibitem[Do~Carmo(2016)]{do2016differential}
Do~Carmo, M.~P.
\newblock \emph{Differential geometry of curves and surfaces: revised and updated second edition}.
\newblock Courier Dover Publications, 2016.

\bibitem[Duran and Grossmann(1986)]{Duran_1986}
Duran, M.~A. and Grossmann, I.~E.
\newblock An outer-approximation algorithm for a class of mixed-integer nonlinear programs.
\newblock \emph{Mathematical Programming}, 36\penalty0 (3):\penalty0 307--339, 1986.

\bibitem[Fosgerau et~al.(2013)Fosgerau, McFadden, and Bierlaire]{FosgMcFaBier13}
Fosgerau, M., McFadden, D., and Bierlaire, M.
\newblock Choice probability generating functions.
\newblock \emph{Journal of Choice Modelling}, 8:\penalty0 1--18, 2013.

\bibitem[Freire et~al.(2016)Freire, Moreno, and Yushimito]{Freire_2016}
Freire, A.~S., Moreno, E., and Yushimito, W.~F.
\newblock A branch-and-bound algorithm for the maximum capture problem with random utilities.
\newblock \emph{European Journal of Operational Research}, 252\penalty0 (1):\penalty0 204--212, 2016.

\bibitem[Gallego and Wang(2014)]{gallego2014multiproduct}
Gallego, G. and Wang, R.
\newblock Multiproduct price optimization and competition under the nested logit model with product-differentiated price sensitivities.
\newblock \emph{Operations Research}, 62\penalty0 (2):\penalty0 450--461, 2014.

\bibitem[Haase(2009)]{Haase_2009}
Haase, K.
\newblock Discrete location planning.
\newblock Technical Report ITLS-WP-09-07, Institute of Transport and Logistics Studies, University of Sydney, 2009.

\bibitem[Haase et~al.(2016)Haase, Müller, Krohn, and Hensher]{Haase_2016}
Haase, K., Müller, S., Krohn, R., and Hensher, D.
\newblock {The maximum capture problem with flexible substitution patterns}.
\newblock Working paper, Hamburg University, 2016.

\bibitem[Haase and Müller(2013)]{Haase_2013}
Haase, K. and Müller, S.
\newblock Management of school locations allowing for free school choice.
\newblock \emph{Omega}, 41\penalty0 (5):\penalty0 847--855, 2013.

\bibitem[Haase and Müller(2014)]{Haase_2014}
Haase, K. and Müller, S.
\newblock A comparison of linear reformulations for multinomial logit choice probabilities in facility location models.
\newblock \emph{European Journal of Operational Research}, 232\penalty0 (3):\penalty0 689--691, 2014.

\bibitem[Haase et~al.(2019)Haase, Kn{\"o}rr, Krohn, M{\"u}ller, and Wagner]{haase2019facility}
Haase, K., Kn{\"o}rr, L., Krohn, R., M{\"u}ller, S., and Wagner, M.
\newblock Facility location in the public sector.
\newblock \emph{Location Science}, pages 745--764, 2019.

\bibitem[Hochbaum and Pathria(1998)]{hochbaum1998analysis}
Hochbaum, D.~S. and Pathria, A.
\newblock Analysis of the greedy approach in problems of maximum k-coverage.
\newblock \emph{Naval Research Logistics (NRL)}, 45\penalty0 (6):\penalty0 615--627, 1998.

\bibitem[Holguin-Veras et~al.(2012)Holguin-Veras, Reilly, and Aros-Vera]{holguin2012new}
Holguin-Veras, J., Reilly, J., and Aros-Vera, F.
\newblock New york city park and ride study.
\newblock Technical report, University Transportation Research Center, 2012.

\bibitem[Kim et~al.(2015)Kim, Pasupathy, and Henderson]{kim2015guide}
Kim, S., Pasupathy, R., and Henderson, S.~G.
\newblock A guide to sample average approximation.
\newblock \emph{Handbook of simulation optimization}, pages 207--243, 2015.

\bibitem[Lamontagne et~al.(2023)Lamontagne, Carvalho, Frejinger, Gendron, Anjos, and Atallah]{lamontagne2023optimising}
Lamontagne, S., Carvalho, M., Frejinger, E., Gendron, B., Anjos, M.~F., and Atallah, R.
\newblock Optimising electric vehicle charging station placement using advanced discrete choice models.
\newblock \emph{INFORMS Journal on Computing}, 35\penalty0 (5):\penalty0 1195--1213, 2023.

\bibitem[Li et~al.(2019)Li, Webster, Mason, and Kempf]{li2019product}
Li, H., Webster, S., Mason, N., and Kempf, K.
\newblock Product-line pricing under discrete mixed multinomial logit demand: winner—2017 m\&som practice-based research competition.
\newblock \emph{Manufacturing \& Service Operations Management}, 21\penalty0 (1):\penalty0 14--28, 2019.

\bibitem[Liu et~al.(2020)Liu, Ma, and Topaloglu]{liu2020assortment}
Liu, N., Ma, Y., and Topaloglu, H.
\newblock Assortment optimization under the multinomial logit model with sequential offerings.
\newblock \emph{INFORMS Journal on Computing}, 32\penalty0 (3):\penalty0 835--853, 2020.

\bibitem[Ljubić and Moreno(2018)]{Ljubic_2018}
Ljubić, I. and Moreno, E.
\newblock Outer approximation and submodular cuts for maximum capture facility location problems with random utilities.
\newblock \emph{European Journal of Operational Research}, 266\penalty0 (1):\penalty0 46--56, 2018.

\bibitem[Mai and Lodi(2020)]{Mai_2020}
Mai, T. and Lodi, A.
\newblock A multicut outer-approximation approach for competitive facility location under random utilities.
\newblock \emph{European Journal of Operational Research}, 284:\penalty0 874--881, 2020.

\bibitem[McFadden and Train(2000)]{McFadden_Train_2000}
McFadden, D. and Train, K.
\newblock Mixed {MNL} models for discrete response.
\newblock \emph{Journal of Applied Econometrics}, 15\penalty0 (5):\penalty0 447--470, 2000.

\bibitem[M{\'e}ndez-Vogel et~al.(2023)M{\'e}ndez-Vogel, Marianov, and L{\"u}er-Villagra]{mendez2023follower}
M{\'e}ndez-Vogel, G., Marianov, V., and L{\"u}er-Villagra, A.
\newblock The follower competitive facility location problem under the nested logit choice rule.
\newblock \emph{European Journal of Operational Research}, 2023.

\bibitem[Miyamoto et~al.(2004)Miyamoto, Vichiensan, Shimomura, and P{\'a}ez]{Miyamoto_2004}
Miyamoto, K., Vichiensan, V., Shimomura, N., and P{\'a}ez, A.
\newblock Discrete choice model with structuralized spatial effects for location analysis.
\newblock \emph{Transportation Research Record}, 1898\penalty0 (1):\penalty0 183--190, 2004.

\bibitem[M{\"u}ller et~al.(2012)M{\"u}ller, Haase, and Seidel]{Muller_2012}
M{\"u}ller, S., Haase, K., and Seidel, F.
\newblock Exposing unobserved spatial similarity: Evidence from german school choice data.
\newblock \emph{Geographical Analysis}, 44\penalty0 (1):\penalty0 65--86, 2012.

\bibitem[Nemhauser and Wolsey(1981)]{nemhauser1981maximizing}
Nemhauser, G.~L. and Wolsey, L.~A.
\newblock Maximizing submodular set functions: formulations and analysis of algorithms.
\newblock In \emph{North-Holland Mathematics Studies}, volume~59, pages 279--301. Elsevier, 1981.

\bibitem[Nemhauser et~al.(1978)Nemhauser, Wolsey, and Fisher]{nemhauser1978analysis}
Nemhauser, G.~L., Wolsey, L.~A., and Fisher, M.~L.
\newblock An analysis of approximations for maximizing submodular set functions—i.
\newblock \emph{Mathematical programming}, 14:\penalty0 265--294, 1978.

\bibitem[Paneque et~al.(2021)Paneque, Bierlaire, Gendron, and Azadeh]{Paneque_2021}
Paneque, M.~P., Bierlaire, M., Gendron, B., and Azadeh, S.~S.
\newblock Integrating advanced discrete choice models in mixed integer linear optimization.
\newblock \emph{Transportation Research Part B: Methodological}, 146:\penalty0 26--49, 2021.

\bibitem[Paneque et~al.(2022)Paneque, Gendron, Azadeh, and Bierlaire]{paneque2022lagrangian}
Paneque, M.~P., Gendron, B., Azadeh, S.~S., and Bierlaire, M.
\newblock A lagrangian decomposition scheme for choice-based optimization.
\newblock \emph{Computers \& Operations Research}, 148:\penalty0 105985, 2022.

\bibitem[Pinsker(1964)]{pinsker1964information}
Pinsker, M.~S.
\newblock \emph{Information and information stability of random variables and processes}.
\newblock Holden-Day, 1964.

\bibitem[Rusmevichientong et~al.(2010)Rusmevichientong, Shen, and Shmoys]{rusmevichientong2010dynamic}
Rusmevichientong, P., Shen, Z.-J.~M., and Shmoys, D.~B.
\newblock Dynamic assortment optimization with a multinomial logit choice model and capacity constraint.
\newblock \emph{Operations research}, 58\penalty0 (6):\penalty0 1666--1680, 2010.

\bibitem[Rusmevichientong et~al.(2014)Rusmevichientong, Shmoys, Tong, and Topaloglu]{rusmevichientong2014assortment}
Rusmevichientong, P., Shmoys, D., Tong, C., and Topaloglu, H.
\newblock Assortment optimization under the multinomial logit model with random choice parameters.
\newblock \emph{Production and Operations Management}, 23\penalty0 (11):\penalty0 2023--2039, 2014.

\bibitem[Salvador and Chan(2004)]{salvador2004determining}
Salvador, S. and Chan, P.
\newblock Determining the number of clusters/segments in hierarchical clustering/segmentation algorithms.
\newblock In \emph{16th IEEE international conference on tools with artificial intelligence}, pages 576--584. IEEE, 2004.

\bibitem[Satopää et~al.(2011)Satopää, Albrecht, Irwin, and Raghavan]{satopaa2011finding}
Satopää, V., Albrecht, J., Irwin, D., and Raghavan, B.
\newblock Finding a ``kneedle'' in a haystack: Detecting knee points in system behavior.
\newblock In \emph{2011 31st international conference on distributed computing systems workshops}, pages 166--171. IEEE, 2011.

\bibitem[Schrijver(2003)]{schrijver2003combinatorial}
Schrijver, A.
\newblock \emph{Combinatorial optimization: polyhedra and efficiency}.
\newblock Springer, 2003.

\bibitem[Vovsha(1997)]{Vovsha_1997}
Vovsha, P.
\newblock \emph{The cross-nested logit model: application to mode choice in the Tel-Aviv metropolitan area}.
\newblock Transportation Research Board, 1997.

\bibitem[Williams(1977)]{Williams_1977}
Williams, H.~C.
\newblock On the formation of travel demand models and economic evaluation measures of user benefit.
\newblock \emph{Environment and planning A}, 9\penalty0 (3):\penalty0 285--344, 1977.

\bibitem[Yaglom and Yaglom(1987)]{gordon1987challenging}
Yaglom, A.~M. and Yaglom, I.~M.
\newblock \emph{Challenging mathematical problems with elementary solutions}, volume~1.
\newblock Courier Corporation, 1987.

\bibitem[Zhang et~al.(2020)Zhang, Rusmevichientong, and Topaloglu]{zhang2020assortment}
Zhang, H., Rusmevichientong, P., and Topaloglu, H.
\newblock Assortment optimization under the paired combinatorial logit model.
\newblock \emph{Operations Research}, 68\penalty0 (3):\penalty0 741--761, 2020.

\bibitem[Zhang et~al.(2012)Zhang, Berman, and Verter]{Zhang_2012}
Zhang, Y., Berman, O., and Verter, V.
\newblock The impact of client choice on preventive healthcare facility network design.
\newblock \emph{OR Spectrum}, 34\penalty0 (2):\penalty0 349–370, 2012.

\end{thebibliography}
\bibliographystyle{plainnat_custom}


\newpage
\appendix
\renewcommand{\thesection}{Appendix \Alph{section}}

\section{Proof of Proposition 2} \label{EC:Prop2Proof}
\proof
The set $X$ of feasible configurations is the same for both models. Therefore, it suffices to show that the restricted models $\widehat{\text{DEQ}}$ and SAA($N,S$) obtained by fixing the decision vector $\boldsymbol{x}$ to a feasible configuration $\bar{\boldsymbol{x}} \in X$ share the same optimal value:
\begingroup
\allowdisplaybreaks
\begin{align}
    \label{model:DEQ_hat_restricted} &\max_{\substack{\boldsymbol{y} \in \{0,1\}^{|\hat{P}|} }}\left\{\sum_{p \in \hat{P}} \hat{\underline{\omega}}_p y_p \Bigg| y_{p} \leq \sum_{d \in D} a^p_{d}\bar{x}_d, \ \forall p \in \hat{P}\right\}, \\
    =&\sum_{p \in \hat{P}} \hat{\underline{\omega}}_p \mathds{1}\left[ \sum_{d \in D} a^p_{d}\bar{x}_d \geq 1 \right], \\
    \label{eq:DEQ_hat_restricted_def_omega} =&\sum_{p \in \hat{P}} \left(\sum_{\substack{(n,s) \in N \times S \\ \boldsymbol{a}^{ns} = \boldsymbol{a}^p}} \frac{q_n}{|S|}\right) \mathds{1}\left[ \sum_{d \in D} a^p_{d}\bar{x}_d \geq 1 \right], \\
    =&\sum_{(n,s) \in N \times S} \frac{q_n}{|S|} \mathds{1}\left[ \sum_{d \in D} a^{ns}_{d}\bar{x}_d \geq 1 \right], \\
    \label{model:SAA_restricted} =&\max_{\substack{\boldsymbol{y} \in \{0,1\}^{|N| \times |S|} }}\left\{\sum_{(n,s) \in N \times S} \frac{q_n}{|S|} y_{ns} \Bigg| y_{ns} \leq \sum_{d \in D} a^{ns}_{d}\bar{x}_d, \ \forall (n,s)\in N \times S \right\}.
\end{align}
\endgroup
The restricted models $\widehat{\text{DEQ}}$ and SAA($N,S$) with $\boldsymbol{x}=\bar{\boldsymbol{x}}$ are respectively given by~(\ref{model:DEQ_hat_restricted}) and~(\ref{model:SAA_restricted}). We obtain (\ref{eq:DEQ_hat_restricted_def_omega}) by replacing the weights $\hat{\underline{\omega}}_p$ with their definition.
\endproof

\section{Illustrative Example of the Knee Detection Method}\label{appendix:knee}
In this section, we illustrate the knee detection method and provide computational insights on the impact of the set partitioning parameter on the computational performance of the \rev{partial Benders decomposition}. We consider randomly generated instances based on the NYC dataset with parameters $\alpha=1$, $\beta=0.1$, and $|S|=1$ scenario.

\begin{figure}[htbp]
\caption{Relative weight of $\hat{P}_1$ and value of $\delta_i$ by fraction of profiles in $\hat{P}_1$}
    \centering
    \includegraphics[width=\textwidth]{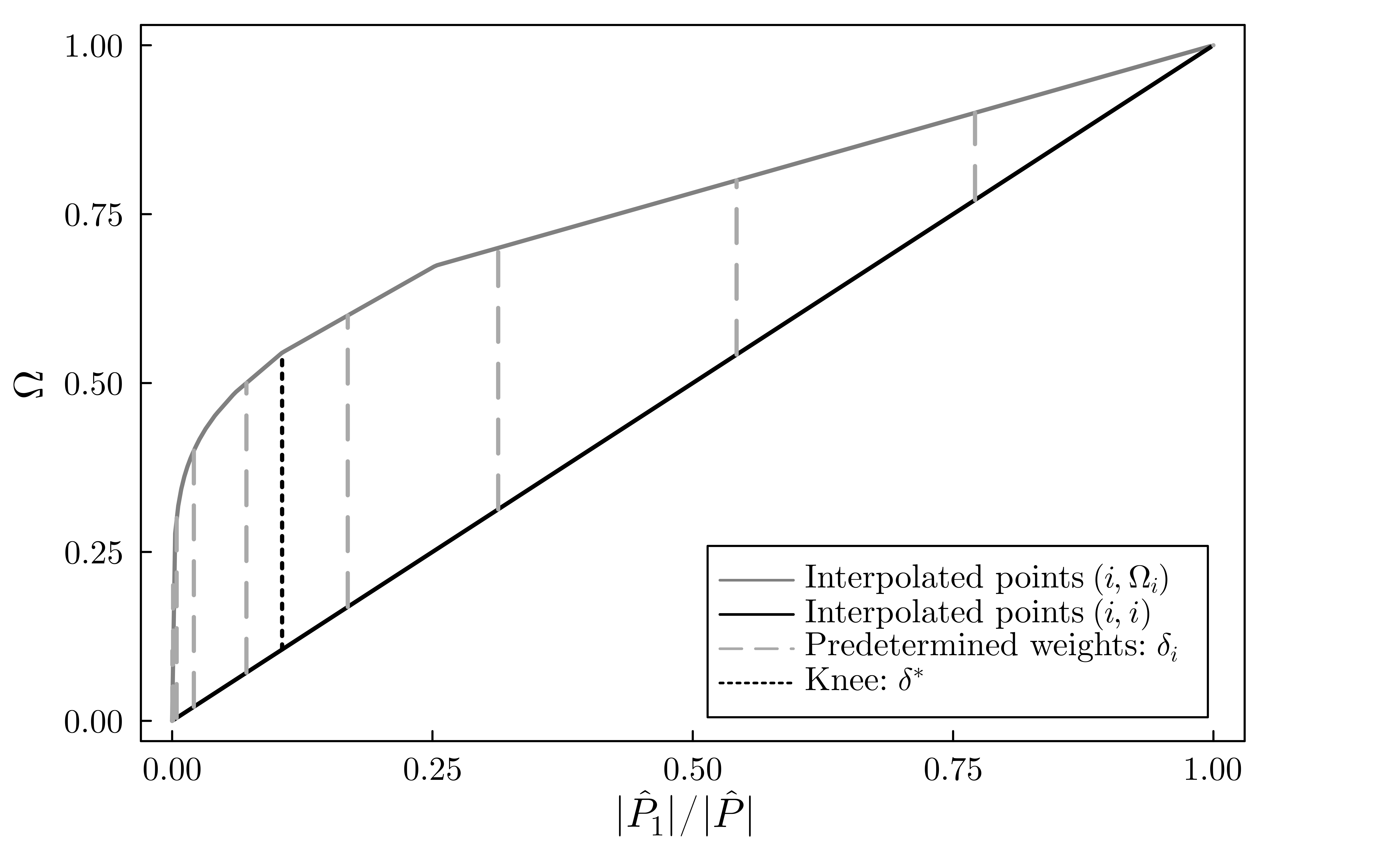}
    \label{fig:knee1}
\end{figure}

Figure~\ref{fig:knee1} compares the points $(i, \Omega_i)$ and the difference $\delta_i = \Omega_i - i/|\hat{P}|$ obtained via the knee detection method and using predetermined weights $\Omega_i \in \{0,0.1,0.2,0.3,0.4,0.5,0.6,0.7,0.8,0.9,1\}$. The knee is reached at point $(i^*, \Omega^*) = (3611, 0.54)$ and achieves a difference $\delta^*=0.44$. In this instance, approximately 10\% of the observed preference profiles thus account for 54\% of the demand. 
\begin{figure}[htbp]
\caption{Number of submodular cuts and $|\hat{P}_1|$ by relative weight of $\hat{P}_1$}
    \centering
    \includegraphics[width=\textwidth]{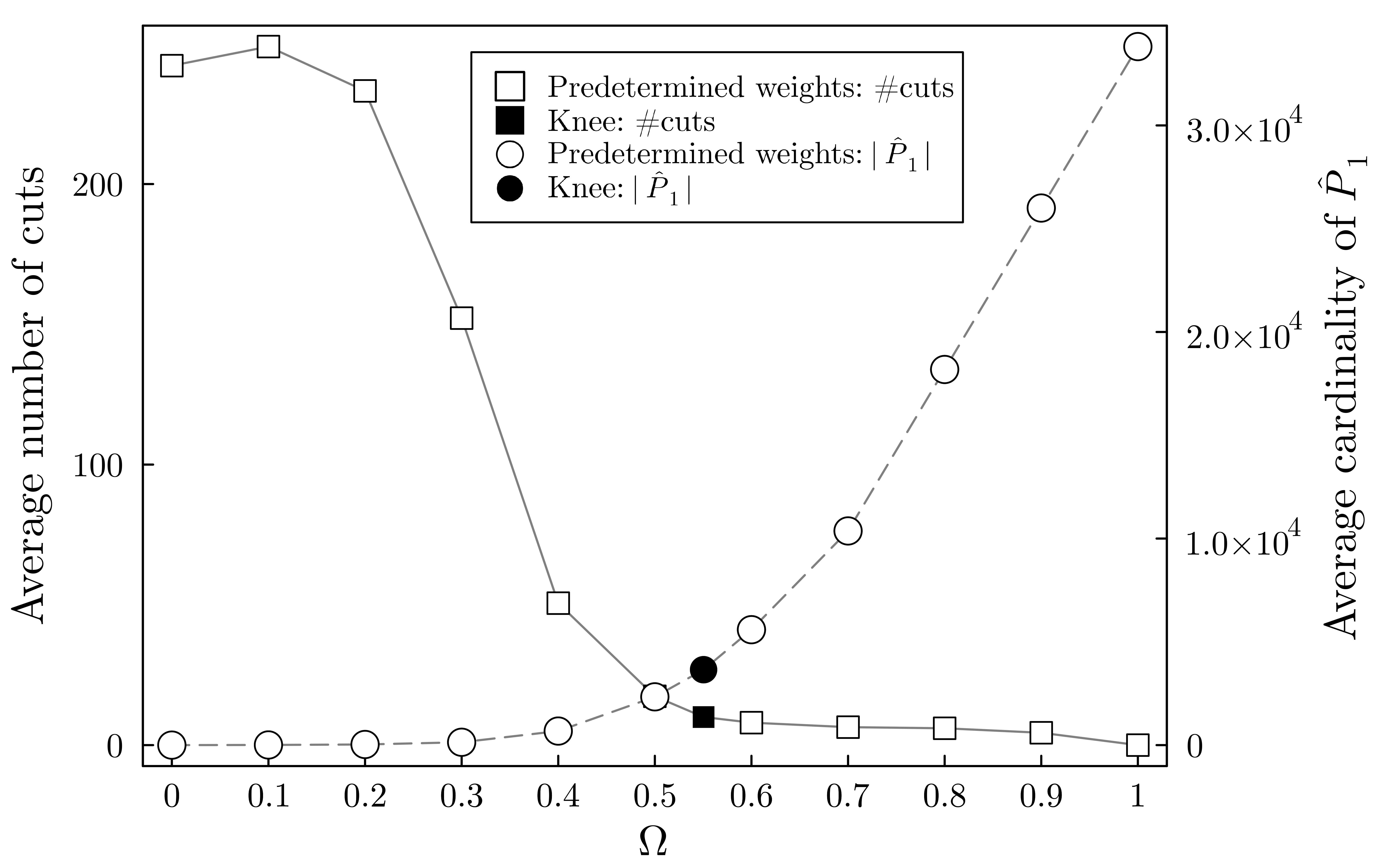}
    \label{fig:knee2}
\end{figure}

\begin{figure}[htbp]
\caption{CPU time by relative weight of $\hat{P}_1$}
    \centering
    \includegraphics[width=\textwidth]{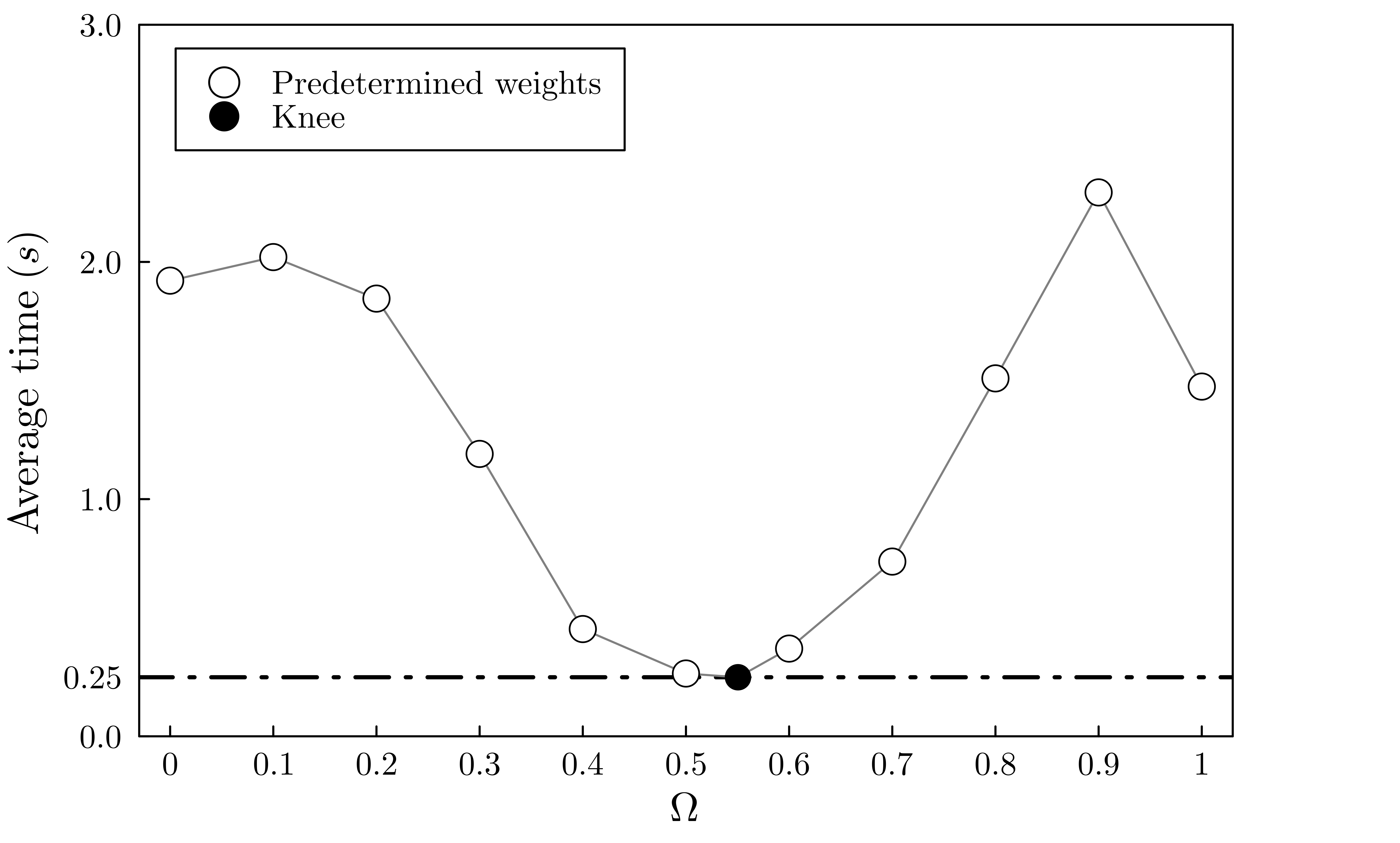}
    \label{fig:knee3}
\end{figure}

Figure~\ref{fig:knee2} depicts the average cardinality of $\hat{P}_1$ and the average number of submodular cuts generated in the B\&C tree when solving model $\widehat{\text{DEQ}}(\hat{P}_1)$ based on the knee detection method and for predetermined weights $\Omega_i$. The average CPU times over 100 instances with a budget of 10 facilities are reported in Figure \ref{fig:knee3}.

It appears that the knee method achieves a good trade-off between the size of the model and the number of generated cuts, with an average of $|\hat{P}_1| = \text{3,651}$ explicitly represented preference profiles and 10 submodular cuts. The pure \rev{Benders decomposition} approach ($\Omega_i=0, \hat{P}_1 = \emptyset$) requires generating 242 cuts on average for these instances. SAAA, which is the other extreme case ($\Omega_i=1, \hat{P}_1 = \hat{P}$), leads to a fairly large 0-1 linear program with $|\hat{P}_1|=\text{33,817}$ demand decision variables on average. These differences in the problem size and the number of generated submodular cuts directly translate into performance differences. Solving $\widehat{\text{DEQ}}$ with \rev{PBD} required 0.25 second on average, compared to 1.72 and 2.74 seconds SAAA and the pure \rev{Benders decomposition} method, respectively.

\section{Entropy and Performance of the \rev{Partial Benders Decomposition Method}}\label{appendix:entropy_performance}
\begin{figure}[htbp]
\caption{Average CPU time for NYC instances of different levels of entropy}
\label{fig:entropy_performance}
\begin{center}
    \includegraphics[width=\textwidth]{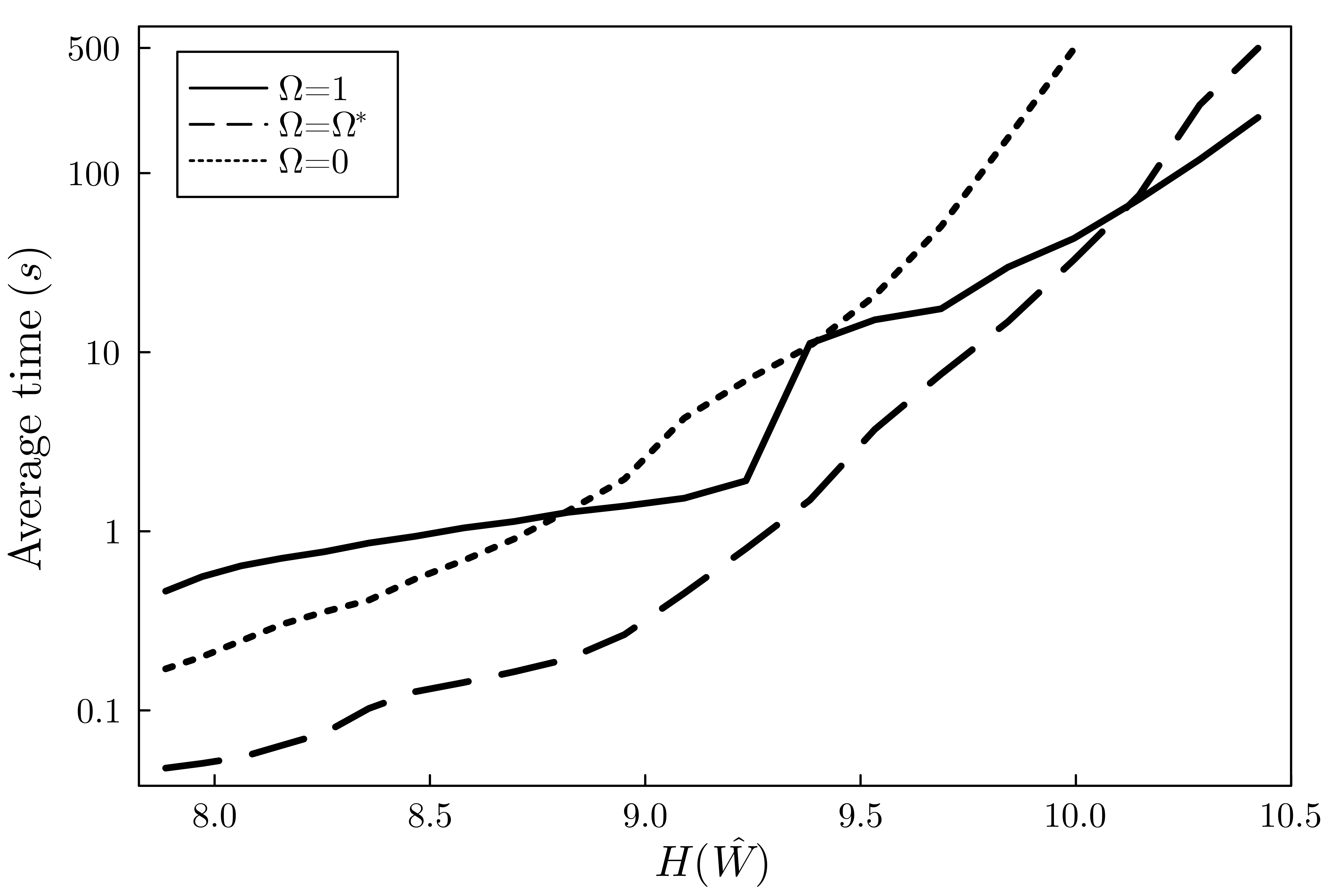}
\end{center}
\end{figure}

Figure \ref{fig:entropy_performance} compares the computational performance of \rev{PBD} ($\Omega=\Omega^*$) with that of SAAA ($\Omega=1$) and the pure \rev{Benders decomposition} method ($\Omega=0$) on instances of different entropy levels. The experimental setup is the same as in \ref{appendix:knee}, except that we consider different values of the parameter $\beta$, which controls the problem's stochasticity level. We solve 10 randomly generated instances with each method for each value of $\beta$ between 0.05 and 0.15, in steps of 0.005. As $\beta$ gets smaller, so does the sensitivity of the customers to distances, which makes their behavior more uncertain (i.e., the relative importance of $\boldsymbol{\varepsilon}$ increases). 

In the lowest entropy setting, \rev{PBD} is the most efficient method, with an average CPU time of 0.05 seconds. The pure \rev{Benders decomposition} method and SAAA take 0.17 and 0.46 seconds per instance on average, respectively. The problem becomes more challenging for all the methods as the entropy increases because so does the number of observed preference profiles, which goes from 19,866 to 60,143 on average between the lowest and highest entropy instances. As expected, the submodularity-based methods are more sensitive to entropy than SAAA. \rev{PBD} ceases to be the most efficient method when $\beta$ exceeds 0.06, and SAAA then provides the best performance. This switch only occurs when the demand is extremely scattered across preference profiles and the empirical distribution of observed preferences approaches the uniform distribution. For $\beta=0.06$, $|\hat{P}|=\text{54,269}$ different nontrivial preference profiles were observed on average, and the resulting empirical distribution $\hat{\boldsymbol{\omega}}$ was associated with an entropy of $H(\hat{W})=10.15$. In comparison, the entropy of the discrete uniform distribution on the same support is 10.90.

\rev{Between SAAA and the pure Benders decomposition, we notice that the latter performs comparatively better when the entropy is quite low.
As shown by \citet{cordeau2019benders}, a pure Benders decomposition method can significantly outperform CPLEX for MCLP instances generated from simple deterministic covering rules. However, as discussed at the end of Section \ref{sec:info_theo:discussion}, such rules lead to a very limited range of preference profiles. In comparison, relying on RUM models (as we do in the context of the CBCFLP) produces a wider variety of preference profiles and can lead to instances with significantly higher entropy depending on the specification of the model and the weight of its random terms. For the NYC instances considered in Figure \ref{fig:entropy_performance}, PBD is approximately one order of magnitude faster than the pure Benders decomposition for every entropy level. On other sets of instances of our testbed, this difference was regularly of two to three orders of magnitude, and the pure Benders decomposition was often outperformed by SAAA. The conclusion is that retaining some information on important profiles in the master problem is critical for the performance of a Benders decomposition approach for the CBCFLP.}

\section{Detailed Presentation of the MIX Dataset} \label{appendix-MIXdataset}
This section provides a complete description of the MIX dataset. In these generative MMNL instances, we assume that there are three types of locations $l \in L$ and three types of customers $k \in \mathcal{K}$. We consider a fixed set of 10 competing facilities (resulting in $|E|=3\cdot 10=30$) and 25 available locations of each type (resulting in $|D|=3\cdot 25=75$), which were generated uniformly on the square $[-15,15]\times[-15,15]$. The distribution of the customers' attributes among the population is given by the random vector $\boldsymbol{\theta} = (\theta_1, \theta_2, K)$, where $(\theta_1, \theta_2)$ denotes a position in the plane, and the categorical random variable $K$ indicates the type of customer. Fixed parameters $\{\delta_{k}\}_{k \in \mathcal{K}}$ and $\{\gamma_{kl}\}_{k \in \mathcal{K}, l \in L}$ respectively control the level of aversion of each customer type to travel distances, and to each type of location. Those are set to $\boldsymbol{\delta} = [3, 1, 2], \boldsymbol{\gamma}_1 = [20, 60, 30], \boldsymbol{\gamma}_2 = [40, 20, 60]$ and $\boldsymbol{\gamma}_3 = [60, 40, 20]$. A random multiplier $\beta \sim U[\beta^-, \beta^+]$ controls the weight of the perceived utility in the function $u_c(\boldsymbol{\theta}, \boldsymbol{\varepsilon}) = v_c(\boldsymbol{\theta}) + \varepsilon_c$. For a facility $c$ of type $l$, the perceived utility is defined as $v_c(\boldsymbol{\theta}) = -\beta(\delta_{K} M_c(\theta_1, \theta_2) + \gamma_{Kl})$, where $M_c(\theta_1, \theta_2)$ denotes the Manhattan distance separating position $(\theta_1, \theta_2)$ from location $c$. The population is distributed across four neighborhoods $j \in \{1,2,3,4\}$. The part of the customers that reside in neighborhood $j$ and the proportion of customers of type $k$ in neighborhood $j$ are specified by two parameters $\pi_j$ and $\rho_{jk}$. The spatial distribution of the population in neighborhood $j$ follows a bivariate normal variable with mean $\boldsymbol{\mu}_j$ and covariance matrix $\Sigma_j$. The following parameters are used in our experiments:
\begin{center}
\vspace{0.5cm}
\scalebox{1}{$
\begin{cases}
\boldsymbol{\pi} &= [0.4, 0.3, 0.2, 0.1] \ ,\vspace{2mm}
\\

[\boldsymbol{\rho}_1, \boldsymbol{\rho}_2, \boldsymbol{\rho}_3, \boldsymbol{\rho}_4] &= \begin{bmatrix} [0.2, 0.7, 0.1], [0.3, 0.4, 0.3], [0.3, 0.4, 0.3], [0.0, 0.2, 0.8]
\end{bmatrix}, \vspace{2mm}
\\

[\boldsymbol{\mu}_1, \boldsymbol{\mu}_2, \boldsymbol{\mu}_3, \boldsymbol{\mu}_4] &= \begin{bmatrix} [2, -2], [-10, -10], [-4, 10], [12, -5]
\end{bmatrix}, \vspace{2mm}
\\

[\Sigma_1, \Sigma_2, \Sigma_3, \Sigma_4] &= \begin{bmatrix}
\begin{bmatrix}
    9\pp & 1\\
    1\pp & 9 \\
\end{bmatrix},
\begin{bmatrix}
    \p 9 \pp & -6\\
    -6 \pp & \p 9 \\
\end{bmatrix},
\begin{bmatrix}
    16\pp& 1\\
    \p 1\pp  & 4 \\
\end{bmatrix},
\begin{bmatrix}
    2\pp & \p 0\\
    0\pp & 21 \\
\end{bmatrix}
\end{bmatrix}.
\end{cases}$}
\vspace{0.8cm}
\end{center}

\begin{figure}[htbp]
\caption{Visualization of the MIX dataset based on a sample of 10,000 customers.}
\label{fig:MMNL-MIX}
\begin{center}
    \includegraphics[width=\textwidth]{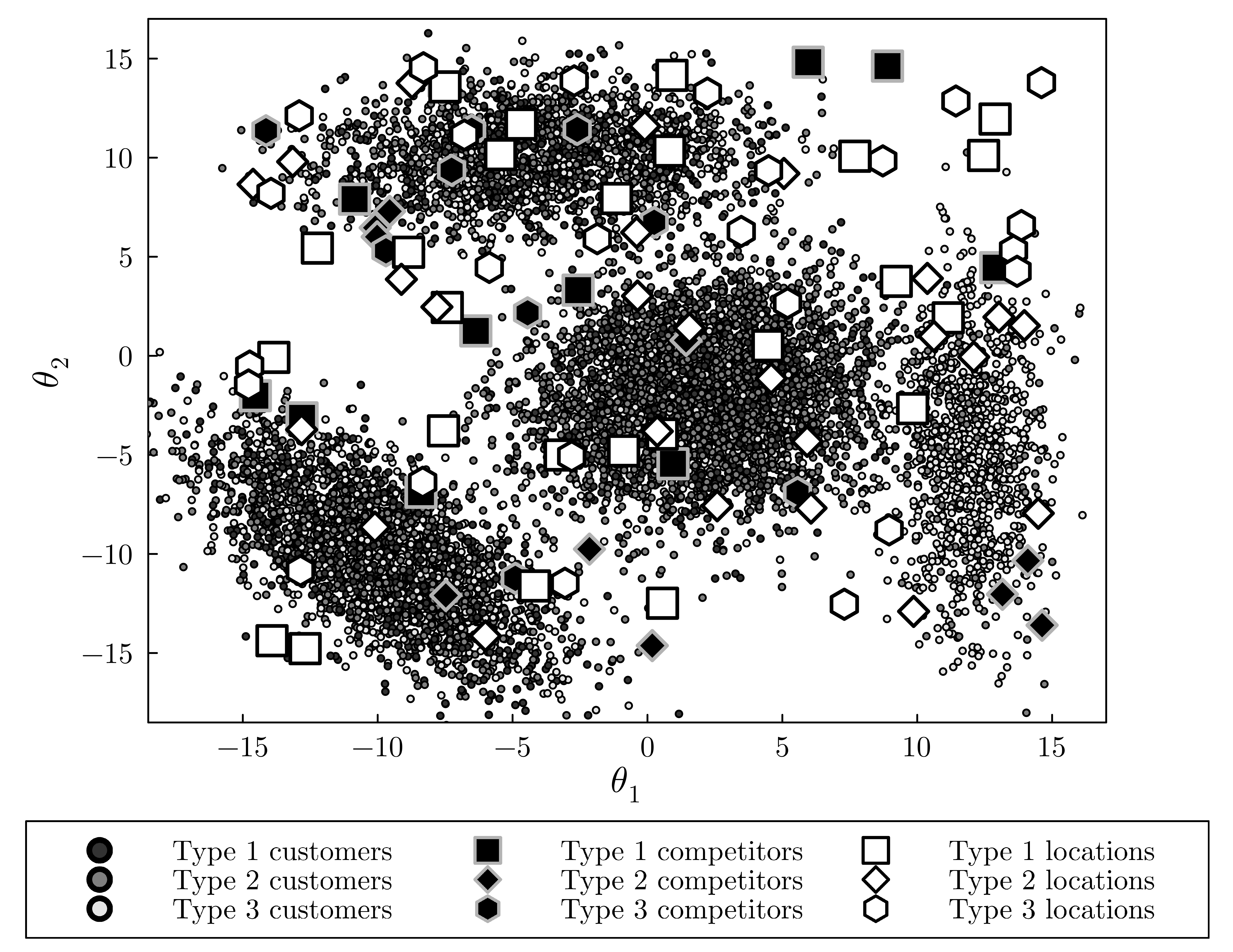}
\end{center}
\end{figure}

Figure~\ref{fig:MMNL-MIX} illustrates the facility locations and the spatial distribution of the customers.

\newpage
\section{Detailed Computational Results.}\label{appendix:results}

This section presents the number of B\&C nodes explored by each simulation-based method and the number of submodular cuts generated by \rev{PBD} in the computational experiments. Tables \ref{table:MNL:BnC} and \ref{table:MMNL:BnC} report the results for the conditional MNL instances and the generative MMNL instances, respectively.

\begin{table*}[htbp]
\centering
\caption{Number of submodular cuts generated and B\&C nodes explored by the simulation-based methods for conditional MNL instances, by entropy level (27 instances per row)}
\label{table:MNL:BnC}
\resizebox{\columnwidth}{!}{
\begin{tabular}{ccccrrrrrrrrrrrr}
\toprule
\multirow{3}{*}{\makecell{Set}} & \multicolumn{3}{c}{\multirow{2}{*}{\makecell{Entropy}}} & \multicolumn{3}{c}{\# submodular cuts} & \multicolumn{9}{c}{\# B\&C nodes} \\
\cmidrule(lr){5-7} \cmidrule(lr){8-16}
& & & &  \multicolumn{3}{c}{\rev{PBD}} &  \multicolumn{3}{c}{\rev{PBD}} &  \multicolumn{3}{c}{SAAA} &  \multicolumn{3}{c}{SAA} \\
\cmidrule(lr){2-4}\cmidrule(lr){5-7}\cmidrule(lr){8-10}\cmidrule(lr){11-13}\cmidrule(lr){14-16} & $\beta$ & $H(\hat{W})$ & $\rev{\log|\hat{P}|}$ &  $|S_1|$  & $|S_2|$ & $|S_3|$ & $|S_1|$  & $|S_2|$ & $|S_3|$ & $|S_1|$  & $|S_2|$ & $|S_3|$ & $|S_1|$  & $|S_2|$ & $|S_3|$ \\
\midrule 
\multirow{4}{*}{\makecell{HM14 \\ $|D|=25$}} 
   & 10 & 3.40 & \rev{\p 4.22} &  8.7 &  6.1 & 7.0 &  4.4 &  3.8 & 4.1 & 0.0 & 0.0 & 0.0 & 0.0 & 0.0 & 0.0 \\
   &  5 & 3.50 & \rev{\p 4.48} &  6.9 & 12.7 & 8.7 &  2.8 &  5.0 & 1.8 & 0.0 & 0.0 & 0.0 & 0.0 & 0.0 & 0.0 \\
   &  2 & 3.87 & \rev{\p 5.25} &  7.2 & 5.5  & 4.1 &  5.2 &  1.5 & 0.4 & 0.0 & 0.0 & 0.0 & 0.0 & 0.0 & 0.0 \\
   &  1 & 4.54 & \rev{\p 6.57} & 20.7 & 11.9 & 7.5 & 37.9 & 11.3 & 5.2 & 0.0 & 0.0 & 0.0 & 0.0 & 0.0 & 0.1 \\
\hline
\multirow{4}{*}{\makecell{HM14 \\ $|D|=50$}} 
   & 10 & 4.40 & \rev{\p 5.10} &  8.3 &  7.0 & 11.7 &  10.1 &  5.5 &  3.7 & 0.0 & 0.0 & 0.0 & 0.0 & 0.0 & 0.0 \\
   &  5 & 4.52 & \rev{\p 5.47} & 11.1 &  6.7 &  7.0 &   9.5 &  2.6 &  1.2 & 0.0 & 0.0 & 0.0 & 0.0 & 0.0 & 0.0 \\
   &  2 & 4.97 & \rev{\p 6.39} & 14.1 & 11.1 &  8.1 &  29.1 & 12.3 &  6.4 & 0.0 & 0.0 & 0.0 & 0.0 & 0.0 & 0.0 \\
   &  1 & 5.86 & \rev{\p 7.74} & 66.1 & 35.0 & 32.0 & 360.0 & 99.7 & 47.6 & 0.0 & 0.0 & 1.8 & 0.0 & 0.0 & 2.6 \\
\hline
\multirow{4}{*}{\makecell{HM14 \\ $|D|=100$}} 
    & 10 & 5.44 & \rev{\p 6.45} &  29.3 &    7.2 &  15.0 &    96.6 &     7.8 &    7.9 &  0.0 &  0.0 &  0.0 &  0.0 &  0.0 &  0.0 \\
    &  5 & 5.76 & \rev{\p 7.04} &  40.1 &   14.4 &   6.9 &   309.9 &    20.6 &    3.6 &  0.0 &  0.0 &  0.0 &  0.0 &  0.0 &  0.0 \\
    &  2 & 6.72 & \rev{\p 8.33} & 249.3 &   93.1 &  64.3 &  7229.4 &   747.1 &  223.7 &  0.0 &  0.0 &  0.0 &  0.0 &  0.0 &  0.0 \\
    &  1 & 8.12 & \rev{\p 9.67} &2990.4 & 1253.6 & 391.0 & 18317.2 & 11281.8 & 3429.6 & 11.8 & 43.5 & 20.1 & 11.3 & 42.8 & 21.1 \\
\hline
\multirow{11}{*}{\makecell{NYC}} & 2.00 & 3.88 & \rev{\p 6.17} & 2.1 &  2.0 & 2.1 &  0.0 & 0.0 & 0.0 & 0.0 & 0.0 & 0.0 & 0.0 & 0.0 & 0.0 \\
                                 & 1.75 & 3.91 & \rev{\p 6.26} &   2.1 &  2.2 & 2.1 &  0.0 & 0.0 & 0.0 & 0.0 & 0.0 & 0.0 & 0.0 & 0.0 & 0.0 \\
                                 & 1.50 & 3.96 & \rev{\p 6.35} &   1.9 &  2.1 & 2.2 &  0.0 & 0.0 & 0.0 & 0.0 & 0.0 & 0.0 & 0.0 & 0.0 & 0.0 \\
                                 & 1.25 & 4.04 & \rev{\p 6.51} &   2.1 &  2.1 & 2.1 &  0.0 & 0.0 & 0.0 & 0.0 & 0.0 & 0.0 & 0.0 & 0.0 & 0.0 \\
                                 & 1.00 & 4.17 & \rev{\p 6.76} &   2.1 &  2.1 & 2.1 &  0.0 & 0.0 & 0.0 & 0.0 & 0.0 & 0.0 & 0.0 & 0.0 & 0.0 \\
                                 & 0.75 & 4.42 & \rev{\p 7.14} &   2.2 &  2.1 & 2.1 &  0.0 & 0.0 & 0.0 & 0.0 & 0.0 & 0.0 & 0.0 & 0.0 & 0.0 \\
                                 & 0.50 & 4.99 & \rev{\p 7.89} &   2.1 &  2.1 & 2.1 &  0.0 & 0.0 & 0.0 & 0.0 & 0.0 & 0.0 & 0.0 & 0.0 & 0.0 \\
                                 & 0.25 & 6.62 & \rev{\p 9.47} &   2.5 &  2.5 & 2.4 &  0.0 & 0.0 & 0.0 & 0.0 & 0.0 & 0.0 & 0.0 & 0.0 & 0.0 \\
                                 & 0.20 & 7.25 & \rev{10.00} &   3.6 &  2.7 & 2.7 &  0.1 & 0.1 & 0.0 & 0.0 & 0.0 & 0.0 & 0.0 & 0.0 & 0.0 \\
                                 & 0.15 & 8.10 & \rev{10.66} &   5.6 &  4.2 & 3.9 &  2.7 & 1.5 & 1.2 & 0.0 & 0.0 & 0.0 & 0.0 & 0.0 & 0.0 \\
                                 & 0.10 & 9.37 & \rev{11.43} &  12.3 & 10.1 & 7.0 & 13.9 & 8.0 & 5.8 & 0.0 & 0.0 & 0.1 & 0.0 & 0.0 & 0.0 \\
\bottomrule
\end{tabular}
}
\end{table*} 

\begin{table*}[htbp]
\centering
\caption{Number of submodular cuts and B\&C nodes generated by the simulation-based methods for the MIX dataset, by entropy level and sample size (25 instances per row)}
\label{table:MMNL:BnC}
\resizebox{\columnwidth}{!}{
\begin{tabular}{ccccccrcc}
\toprule
\multicolumn{4}{c}{Entropy} & \multirow{2}{*}{\makecell{$|N|$\\(thousands)}} & \# submodular cuts & \multicolumn{3}{c}{\# B\&C nodes} \\
\cmidrule(lr){1-4}\cmidrule(lr){6-6}\cmidrule(lr){7-9} $\beta^-$ & $\beta^+$ & $H(\hat{W})$ & $\rev{\log|\hat{P}|}$ &  &  \rev{PBD} & \rev{PBD} & SAAA & SAA  \\
\midrule 
\multirow{5}{*}{1} & \multirow{5}{*}{2} 
   & 5.09 & \rev{\p 7.14} &   16 &   15.1  & 3.2  & 0.0  & 0.0 \\
 & & 5.11 & \rev{\p 7.36} &   32 &    \p8.8  & 1.6  & 0.0  & 0.0 \\
 & & 5.12 & \rev{\p 7.56} &   64 &    \p3.7  & 0.0  & 0.0  & 0.0 \\
 & & 5.12 & \rev{\p 7.77} &  128 &   10.3  & 3.0  & 0.0  & 0.0 \\
 & & 5.12 & \rev{\p 7.95} &  256 &   15.8  & 3.2  & 0.0  & 0.0 \\
\hline\multirow{5}{*}{0.5} & \multirow{5}{*}{1} 
   & 5.61 & \rev{\p 7.91} &  16 &  10.2  & 7.2  & 0.0  & 0.0 \\
 & & 5.63 & \rev{\p 8.23} &  32 &   \p8.0  & 3.4  & 0.0  & 0.0 \\
 & & 5.65 & \rev{\p 8.54} &  64 &   \p4.9  & 1.2  & 0.0  & 0.0 \\
 & & 5.65 & \rev{\p 8.83} & 128 &   \p3.8  & 0.4  & 0.0  & 0.0 \\
 & & 5.67 & \rev{\p 9.11} & 256 &  10.1  & 1.5  & 0.0  & 0.0 \\
\hline\multirow{5}{*}{0.25} & \multirow{5}{*}{0.5}  
   & 6.63 & \rev{\p 8.93} &  16 &  13.3  & 12.5  & 0.0  & 0.0 \\
 & & 6.68 & \rev{\p 9.36} &  32 &  11.2  & 11.6  & 0.0  & 0.0 \\
 & & 6.74 & \rev{\p 9.79} &  64 &   \p8.7  &  7.6  & 0.0  & 0.0 \\
 & & 6.77 & \rev{  10.19} & 128 &   \p6.5  &  3.7  & 0.0  & 0.0 \\
 & & 6.80 & \rev{  10.59} & 256 &   \p5.7  &  2.6  & 0.0  & 0.0 \\
\hline\multirow{5}{*}{0.125} & \multirow{5}{*}{0.25} 
   & 7.86 & \rev{\p 9.77} &  16 &  47.5  & 166.0  & 0.0  & 0.0 \\
 & & 8.01 & \rev{  10.30} &  32 &  40.2  &  94.9  & 0.3  & 0.0 \\
 & & 8.16 & \rev{  10.83} &  64 &  22.9  &  48.9  & 0.2  & 1.8 \\
 & & 8.29 & \rev{  11.36} & 128 &  23.8  &  33.7  & 1.0  & 1.4 \\
 & & 8.39 & \rev{  11.89} & 256 &  23.6  &  26.9  & 0.0  & 0.0 \\
\bottomrule
\end{tabular}
}
\end{table*}

Two main observations can be drawn from these tables. First, both SAAA and SAA explore a negligible number of nodes, meaning that practically all the solving time is spent at the root node for these methods. Second, the number of B\&C nodes and the number of submodular cuts generated by \rev{PBD} generally decrease with the number of simulated customers. For example, as reported in Table \ref{table:MNL:BnC}, the average number of submodular cuts for the HM14 instances with $|D|=100$ decreases from 2,990.4 to 391.0 and the number of B\&C nodes decreases from 18,317.2 to 3,429.6 between $|S_1|=10$ and $|S_3|=\text{1,000}$ scenarios. The results of Table \ref{table:MMNL:BnC} indicate the same trend with respect to the number $|N|$ of customers. 

\end{document}